\documentclass[12pt]{article}

\newcommand{\comment}[1]{}
 
\usepackage{amssymb}
\usepackage{amsthm}
\usepackage{amsmath}

\usepackage{tikz}
\usetikzlibrary{decorations.markings}

\theoremstyle{plain}
\newtheorem{theorem}{Theorem}[section]
\newtheorem{lemma}[theorem]{Lemma}
\newtheorem{corollary}[theorem]{Corollary}

\theoremstyle{definition}
\newtheorem{example}[theorem]{Example}
\newtheorem{conjecture}[theorem]{Conjecture}
\newtheorem{question}[theorem]{Question}

\newcommand{\cE}{{\cal E}}
\newcommand{\cB}{{\cal B}}
\newcommand{\cG}{{\cal G}}
\newcommand{\cS}{{\cal S}}
\newcommand{\cO}{{\cal O}}
\newcommand{\cX}{{\cal X}}
\newcommand{\ttt}{Tic-Tac-Toe}
\newcommand{\MB}{Maker-Breaker}
\newcommand{\entryi}{\alpha}
\newcommand{\entryii}{\beta}

\newcommand{\TD}{{\rm TD}}
\newcommand{\BIBD}{{\rm BIBD}}

\newcommand{\disunion}{\, \dot\cup \,}

\title{\ttt{} on Designs}
\author{
Peter Danziger\thanks{Department of Mathematics, 
	Toronto Metropolitan University, Ontario, Canada}
 \thanks{Supported by NSERC Discovery Grant RGPIN-2027-03816}
\and
Melissa A. Huggan\thanks{Department of Mathematics, Vancouver Island University, Nanaimo, British Columbia, Canada}
\thanks{This author was partially funded by the Atlantic Association for Research in the Mathematical Sciences (AARMS) and Mount Allison University.} 
\and
Rehan Malik$^*$
\and
Trent G.\ Marbach$^*$
}

\begin{document}
\maketitle

\begin{abstract}
	We consider playing the game of \ttt{} on block designs BIBD$(v,k,\lambda)$ and transversal designs TD$(k, n)$.
	Players take turns choosing points and the first player to complete a block wins the game. We show that triple systems, BIBD$(v,3,\lambda)$, are a first player win if and only if $v \geq 5$. Further, we show that for $k=2,3$, TD$(k, n)$ is a first player win if and only if $n\geq k$.
	We also consider a {\em weak} version of the game, called \MB{}, in which the second player wins if they can stop the first player from winning.
	In this case, we adapt known bounds for when either the first or second player can win on BIBD$(v,k,1)$ and TD$(k,n)$, and show that for \MB{}, BIBD$(v,4,1)$ is a first player win if and only if $v \geq 16$.
	We show that TD$(4,4)$ is a second player win, and so the second player can force a draw in the regular game by playing the same strategy.
\end{abstract}

\noindent
{\bf Keywords:} Positional Games, \ttt{}, \MB{}, Games on Block Designs, Games on  Transversal designs, Games on Triple Systems.

\section{Introduction} 
\label{sec:intro}

In its most general setting, \ttt{} is a positional game played by two players on a hypergraph $(V, {\cE})$, where $V$ is a set of \emph{vertices} and $\cE$ is a set of subsets of $V$ called hyperedges.
The players take turns choosing vertices that have not been chosen yet, and the first player to occupy all of the vertices of a hyperedge is the winner.
The convention has arisen that the first player is referred to as Xeno and the second player as Ophelia.
Positional games were first introduced by Erd\H{o}s and Selfridge~\cite{ES} and were further considered by Berge \cite{Berge} and Beck~\cite{Beck81}.
 Erd\H{o}s and Selfridge~\cite{ES} noted that by strategy stealing, it is impossible for Ophelia to win the game, so the best outcome she can hope for is a draw.
A solution to the game for Xeno is a strategy that allows him to win the game, and a solution for Ophelia is a strategy that allows her to always force a draw.
We say that a hypergraph is {\em Xeno win} if Xeno has a winning strategy, regardless of Ophelia's strategy. 
A hypergraph is called {\em Ophelia draw} if it is not Xeno win. 
\ttt{} has been explored in other structures such as graphs \cite{Beeler18} and hypercubes \cite{GH02}, as well as affine and projective planes \cite{CarrollDougherty}.  For an overview of the subject, see \cite{Beck05, Beck08}.

\subsection{Hypergraphs and Designs}
\label{subs Designs}

We say that a hypergraph is \emph{$k$-uniform} if all hyperedges are of the same cardinality, $k$, and \emph{linear} if no pair of vertices from $V$ appears in more than one hyperedge. 
Linear hypergraphs are closely related to the study of \emph{designs}, and this paper aims to study \ttt{} on designs and design-like objects. 

A \emph{design} is a pair $(X, \cB)$, where $X$ is a set of \emph{points} and $\cB$ is a set of subsets of $X$, called \emph{blocks}, such that every pair of points of $X$ appears exactly $\lambda$ times in the blocks of $\cB$. 
A \emph{Balanced Incomplete Block Design} BIBD$(v, k, \lambda)$ is a design where there are $v$ points, so $|X|=v$; each block has size $k$, so $|B|= k$ for every $B\in\cB$, and each pair of points occur in exactly $\lambda$ blocks. 
Particular sub-classes of BIBDs that we will focus on are \emph{Triple Systems} TS$(v,\lambda)$, which are BIBD$(v,3,\lambda)$, and \emph{Steiner Triple Systems} STS$(v)$, which are BIBD$(v,3,1)$. 
A design can be thought of as a hypergraph, where points correspond to vertices and blocks correspond to hyperedges. 

It is well known that necessary conditions for the existence of a BIBD$(v,k,\lambda)$ \cite{Handbook} are that 
\begin{equation}
\label{bibd nec}
\lambda (v-1) \equiv 0 \bmod (k-1)
\mbox{ and } \lambda v(v-1)\equiv 0\bmod k(k-1).
\end{equation}
When $\lambda=1$, this can be summarized as $v\equiv 1$ or $k\bmod{k(k-1)}$. Values of $(v,k,\lambda)$ which satisfy Equation~(\ref{bibd nec}) are called \emph{admissible}.

We note that in a BIBD$(v,k,\lambda)$, there are $|X|=v$ points, $|\mathcal{B}|=\frac{\lambda v(v-1)}{k(k-1)}$ blocks, each pair of points occur in $\lambda$ blocks, and each point is in $r=\frac{\lambda(v-1)}{k-1}$ blocks, $r$ is known as the replication number. 
Note that when playing \ttt{} on a BIBD$(v,k,\lambda)$, as with any design, the blocks are winning sets. 

A \emph{Transversal Design} \TD$(k,n)$ is a triple $(X, \cG, \cB)$ where: $X$ is a set of $kn$ points; $\cG$ is a partition of $X$ into $k$ subsets of cardinality $n$, called \emph{groups}; and $\cB$ is a set of subsets of $X$, each of cardinality $k$, called \emph{blocks}, such that every pair of points not in the same group appears exactly once in some block. We note that when playing \ttt{} on a transerval design, the blocks are winning sets, but the groups are not.

A \emph{parallel class} of a design is a set of blocks that partition the design's point set.
A design is called \emph{resolvable} if the block set can be partitioned into parallel classes. 
A \emph{Resolvable Transversal Design} RTD$(k,n)$ is a resolvable \TD$(k,n)$. 
The block set of the \emph{affine plane} of order $n$, written $\pi_n$, can be obtained from an RTD$(n,n)$ by defining the set of lines of the plane to be the union of the transversal design's groups and blocks.
Further, a \emph{projective plane} of order $n$, denoted $\Pi_n$, can be obtained by adding an extra point $u$ to a \TD$(n+1,n)$, and adding the block $g\cup \{u\}$ for each group $g$ of the design. 
It is also known that an RTD$(k,n)$ can be used to construct a TD$(k+1,n)$ by connecting a new vertex to each block in a parallel class, for each parallel class in the RTD$(k,n)$.
These processes are reversible, and so the existence of an RTD$(n,n)$,  a \TD$(n+1,n)$, a $\pi_n$, and a $\Pi_n$, are all equivalent.

It is well known that if a \TD$(k,n)$ exists, then $k \leq n+1$.
While the existence of a TD$(k,n)$, $k\leq n+1$, is known for all prime power orders of $n$, there are many pairs $(k, n)$ with $k\leq n+1$ for which the existence of a TD$(k, n)$ is unknown. Henceforth, when dealing with the general case, we will always assume that the required design exists, without further mention.
We refer the reader to \cite{Handbook} and \cite{Stinson} for further information, results, and terminology on designs.

\subsection{\MB{}: A weak variant of \ttt{}}

It has been noted that solving the game of \ttt{} is hard in general, see \cite{Beck05,Beck08}. 
As a result of the difficulty of finding solutions to even relatively simple situations, a weak version of the game is sometimes considered, which is known as the \emph{\MB{}} game. 
This has also been referred to as the positional game of the first type by Berge \cite{Berge}, who also referred to \ttt{} as the positional game of the second type. 
We say that a hypergraph is {\em Maker win} if Maker has a winning strategy, regardless of Breaker's strategy and similarly for {\em Breaker win}. 

In the \MB{} version of the game, the two players take the roles of \emph{Maker} and \emph{Breaker}, which correspond to Xeno and Ophelia, respectively.
Maker plays first and wins the game by completing a hyperedge, as usual.
Breaker wins the game if all vertices have been chosen, but Maker has not won. 
Thus, Breaker wins by stopping Maker from completing a winning hyperedge, but cannot win the game by filling a winning hyperedge themselves. 
It is worth noting that the introduction of \MB{} fundamentally changes the nature of the game as Maker does not have to worry about potential threats from Breaker.

We note that if Breaker can win in the \MB{} version, then Ophelia can force a draw in the \ttt{} version by playing the same strategy. 
On the other hand, if Xeno can win in the \ttt{} version, then Maker can win in the \MB{} version by playing the same strategy.
We thus have the following theorem.

\begin{theorem}
    \label{MB to ttt}
    If Breaker can win \MB{} on a given hypergraph, then Ophelia can force a draw when playing \ttt{} on the same hypergraph.
    Similarly, if Xeno can win \ttt{} on a given hypergraph, then Maker can win \MB{} on the same hypergraph.
    \end{theorem}
We note that the converse of Theorem~\ref{MB to ttt} is in general false. 
That is, a Maker win in \MB{} does not imply Xeno can win when playing \ttt{} on the same hypergraph, as Ophelia may be able to force Xeno into making certain moves in \ttt{}. 
Similarly, if Ophelia can force a draw in \ttt{}, she may have used the ability to force moves to achieve this, and so Breaker is not guaranteed to be able to win \MB{} on the same hypergraph. 
Indeed, we have verified through a computer search that a TD$(4,5)$ is an Ophelia draw but is a Maker win. As far as we know, this is the first known example of such a case.

Nonetheless, some of the strategies for the game of \ttt{} can be used to prove an equivalent result for the \MB{} game.
In particular, the arguments involving weight functions and scores described below hold for both games.

\subsection{Previous results}

\subsubsection*{Affine and projective plane results}
Carroll and Dougherty \cite{CarrollDougherty} considered \ttt{}  on affine and projective planes, we summarize their results here. 
\begin{theorem}[\cite{CarrollDougherty}]
\label{Affine and Projective Planes Result}
Ophelia can force a draw when playing \ttt{} on: 
\begin{itemize}
    \item an affine plane of order $n$, $\pi_n$, if and only if $n \geq 5$;
    \item a projective plane of order $n$, $\Pi_n$, if and only if $n \geq 3$.
\end{itemize}
\end{theorem}

While the authors do not explicitly consider the \MB{} variant, they do use weight functions (see below) and many of their proofs are applicable to the \MB{} game. 
One exception to this is the case where \MB{} is played on a $\pi_5$, which cannot be derived from their work on \ttt{}, as they relied on a computer search to check that $\pi_5$ is an Ophelia draw. We have performed this computation for \MB{} and found that Breaker can also win on $\pi_5$. 
We thus have the following result.

\begin{theorem}
\label{Affine and Projective Planes Result - Maker Breaker}
Breaker can win when playing \MB{} on: 

\begin{itemize}
    \item an affine plane of order $n$, $\pi_n$, if and only if $n\geq 5$;
    
	\item a projective plane of order $n$, $\Pi_n$, if and only if $n \geq 3$.
 \end{itemize}
\end{theorem}

\subsubsection*{General hypergraph results}

While we will primarily be studying block designs in this work, 
we may think of the block design as a hypergraph, where blocks correspond to hyperedges and points to vertices. 
Working with hypergraphs in general can provide us with a flexible tool that yields results for designs. 
We will use the terms hyperedge and vertex when explicitly working on hypergraphs, as opposed to using the terms block and point when working on designs.

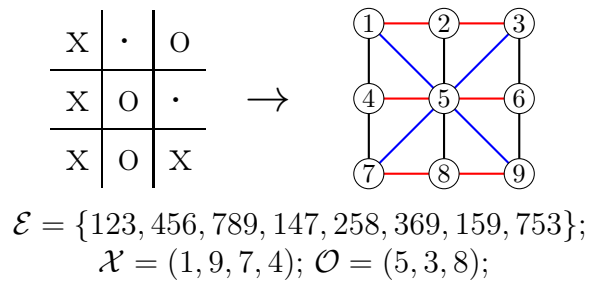
\begin{figure}
\begin{center}
	{\Large
\begin{tabular}{p{1ex}|p{1ex}|p{1ex}}
        x & $\cdot$ & o \\\hline x&o&$\cdot$ \\\hline x&o&x \\
\end{tabular}
\hspace{10pt}$\rightarrow$\hspace{10pt}
}
\raisebox{-6ex}{
  \begin{tikzpicture}[scale=1]
    \draw[thick] (0,0) -- (0,2);
    \draw[thick] (1,0) -- (1,2);
    \draw[thick] (2,0) -- (2,2);
    \draw[red,thick]  (0,0) -- (2,0);
    \draw[red,thick] (0,1) -- (2,1);
    \draw[red,thick] (0,2) -- (2,2);
    \draw[blue,thick] (0,0) -- (2,2);
    \draw[blue,thick] (0,2) -- (2,0);
    
    \foreach \x in {0,1,2}
    \foreach \y in {0,1,2}
    {
    \draw[fill=white] (\x,\y) circle (0.2cm);
    }
    \fill (0,0) node {\footnotesize 7};
    \fill (0,1) node {\footnotesize 4};
    \fill (0,2) node {\footnotesize 1};
    \fill (1,0) node {\footnotesize 8};
    \fill (1,1) node {\footnotesize 5};
    \fill (1,2) node {\footnotesize 2};
    \fill (2,0) node {\footnotesize 9};
    \fill (2,1) node {\footnotesize 6};
    \fill (2,2) node {\footnotesize 3};
  \end{tikzpicture}
}

\ 

  $\cE=\{123,456,789,147,258,369,159,753\}$;
  
  $\cX=(1,9,7,4)$;   
  $\cO=(5,3,8)$;
  \vspace*{-2ex}
\end{center}
\caption{A game of \ttt{} played on a standard $3\times 3$ board. The winning hyperedges and the state of a proposed game
 are also given.}
\label{fig:tictactoe}
\end{figure}
When playing \ttt{} on some hypergraph, we define the \emph{state} of the game to be the pair of move sets, $\cS = (\cX , \cO)$, where $\cX = (X_1,\ldots, X_r)$ is the ordered set of vertices played by Xeno and
$\cO = (O_1, \ldots, O_s)$ is the ordered set of vertices played by Ophelia. 
Note that necessarily $s \in \{r - 1, r\}$. 
For such a state, we say that there have been $r + s$ \emph{moves} and that we are in the $r^{\rm th}$ \emph{round}. 
Ophelia will take the next move if $s = r - 1$, and Xeno otherwise. 
See Figure~\ref{fig:tictactoe} for an example of a state of a \ttt{}  game played on a standard $3\times 3$ board. 
The \emph{state} of the \MB{} game is defined analogously, with $\cX = (X_1,\ldots, X_r)$ now being the ordered set of vertices played by Maker and
$\cO = (O_1, \ldots, O_s)$ the ordered set of vertices played by Breaker. 

Given a state $(\cX , \cO)$ of either game, Erd\H{o}s and Selfridge \cite{ES} introduced the notion of a weighting function for each hyperedge $h$:
\[
	w(h) = \left\{
	\begin{array}{ll}
		2^{|\mathcal{X} \cap h| - |h|} &\text{\rm if $h$ is disjoint from } \mathcal{O} \text{\rm , and } \\
		0 &\text{\rm otherwise.}
	\end{array}
	\right.
\]
We can define the weight of a vertex $u\in V$ to be the sum of the weights of the hyperedges in which that vertex appears, $w(u) = \sum_{B: u\in B} w(B)$.
Finally, we define the \emph{score} of a state to be the sum of the weights of all the hyperedges, $\sum_{h\in \cE} w(h)$. 

A few properties of this scoring function can help to understand the decision behind its choice. 
Firstly, when playing \MB{} (\ttt{}), any hyperedge that Breaker (Ophelia) has played in cannot be a winning hyperedge for Maker (Xeno), and hence contributes $0$ to the score. 
Secondly, any hyperedge of current score $s$ will be worth $s+s$ if Maker (Xeno) plays in it, and worth $s-s$ if Breaker plays in it. 
Thus the weight of every hyperedge is either unchanged or doubles in value when Maker (Xeno) plays, and the weight of every hyperedge is either unchanged or goes to zero when Breaker (Ophelia) plays. 
When Breaker (Ophelia) takes her turn in order to minimise the score, say decreasing the score by $s$, then in Maker's (Xeno's) next turn, the score can increase by at most $s$. 
Thus, in this case, the score at Breaker's (Ophelia's) next move is at most the score at her current one.
Finally, if Maker (Xeno) wins on hyperedge $B$, then in the winning state, $w(B)=1$, so the score of a winning state is at least $1$. 
These observations lead to the following theorem.
\begin{theorem}[\cite{Beck81, ES}]
\label{score<1}
If the score of a state in \MB{} (\ttt{})  after Maker's (Xeno's) turn is ever below $1$, then Breaker is able to win (Ophelia can force a draw). 
\end{theorem}

We can then categorize those hypergraphs for which during a \MB{} (\ttt{}) game, Maker's score after its first turn is below 1, meaning that Breaker can win (Ophelia can force a draw) on these hypergraphs. 
\begin{theorem}[\cite{Beck81}]\label{thm:MB_general_breaker}
	If $\mathcal{H}$ is a $k$-uniform linear hypergraph that satisfies 
	\[
    	2^ k > |\mathcal{E}| + \max_{u\in V} |\{h \in \mathcal{E} : u \in h\}|, 
    \]
    then Breaker can win in \MB{} (Ophelia can force a draw in \ttt{}) on $\mathcal{H}$.
\end{theorem}
\begin{proof}
    The score at the start of play is $2^{-k} |\cE|$. When Xeno plays his first move on some vertex, say $u$, each block containing $u$ will double its weight, going from weight $2^{-k}$ to weight $2^{-k+1}$. Thus the score after this move is at most $2^{-k}(|\mathcal{E}| + \max_{u\in V} |\{h \in \mathcal{E} : u \in h\}|)$, and Theorem~\ref{score<1} finishes the proof. 
    The result for \ttt{} follows from Theorem~\ref{MB to ttt}.
\end{proof}

We note that in the case of a design, the term $\max_{u\in V} |\{h \in \mathcal{E} : u \in h\}|$ is the replication number, $r$, of a point.
A companion result to Theorem~\ref{thm:MB_general_breaker} for \MB{} in the Maker win cases was proved by Beck~\cite{Beck81}, adapted here for uniform linear hypergraphs.
\begin{theorem}[\cite{Beck81}]
\label{MB win condition}
Suppose that the \MB{} 
game is played on a $k$-uniform linear hypergraph $(V, \cE)$.  Then Maker can win whenever
\[
2^{k-3}|V| < |\cE|.
\]
\end{theorem}

\subsection{Optimal strategies}

The discussion of weight and score above suggests a \emph{score-optimizing} strategy that either player can adopt, where the player always chooses a vertex of maximum weight. However, in general, only using score-optimizing as a strategy can fail to block potential wins and can fail to make a winning move that is available. 
Therefore, when referring to a \emph{score-optimizing strategy}, we mean a strategy that prioritizes playing a winning move, and if none exist, prioritizes blocking any immediate wins of their opponent, and if there are no such wins, then chooses a vertex of maximum weight, breaking ties by any method. In previous works, there has been a heavy reliance on using score-optimizing strategies, under an implicit assumption that these are either optimal or near-optimal strategies. We thus ask the following question.
\begin{question}\label{qn:score_opt}
In either the \ttt{} or \MB{} games, is a score-optimizing strategy always an optimal strategy for either player?
\end{question}

We will show that this is false, at least for the \ttt{} game. 
There are hypergraphs $\mathcal{H}$ where one player could win on $\mathcal{H}$ by playing an optimal strategy, but will lose if they, or even both players, play a score-optimizing strategy. 
For a simple example, consider playing \ttt{} on the standard $3\times 3$ board shown in Figure \ref{fig:tictactoe}, with state $\mathcal{X} = (1,9,7,4)$ and $\mathcal{O} = (5,3,8)$ as shown. 
Xeno wins with hyperedge $\{1, 4, 7\}$ and it is easy to verify that Ophelia has played a score-optimizing strategy by making plays that minimize the score on each move. 
However, it is well known that Ophelia can force a draw on the standard \ttt{} game. 
We note that in this example, Xeno does not play using a score-optimizing strategy. 

For designs in particular, we show in Theorem~\ref{k=3 is Xeno win} that any \TD$(3,n)$ is a Xeno win.  
In contrast to this, Theorem~\ref{pan-hamiltonian} shows that Ophelia can force a draw on certain transversal designs if Xeno plays a score-optimizing strategy (in fact, Ophelia also plays a score-optimizing strategy in this case).
This suggests an alternate question. 

\begin{question}
For which designs and hypergraphs does there exist a score-optimizing strategy that is optimal?
\end{question}

\subsection{Summary}
Previous works have provided a number of general results such as Theorem~\ref{thm:MB_general_breaker} and Theorem~\ref{MB win condition} for determining which player will win \MB{} or \ttt{} under certain conditions. 
It has been noted that identifying the winning player in either game outside of these cases is, in general, difficult to do. 
In this paper, we aim to identify some classes of designs which are not amenable to these results, and so may be described as tough cases to analyse. 
In particular, we consider BIBDs and TDs with small block size. 
Then, to address these cases, we develop techniques to solve the problem on some of these designs. 
Additionally, we shall answer Question~\ref{qn:score_opt} when \ttt{} is played on designs.  

In Section~\ref{Weak Breaker win}, we consider the \MB{} variant of the game. 
We focus on \MB{} games in which Breaker can win.
These will be the most useful to us, as a Breaker winning strategy in the \MB{} game can be converted by Theorem~\ref{MB to ttt} to a strategy where Ophelia can force a draw in \ttt{} on the same structure. 
We classify how Theorem~\ref{thm:MB_general_breaker} and Theorem~\ref{MB win condition} apply to BIBDs in Subsection~\ref{subs:MB_BIBD}, and transversal designs in Subsection~\ref{subs:MB_TD}. 
In the process we completely classify \MB{} on BIBD$(v,4,1)$ in Theorem~\ref{BIBD(v,4,1)}.
A standout result of Section~\ref{Weak Breaker win} is a proof that Breaker can win on a \TD$(4,4)$, the smallest case where these results do not apply.  Note that this means that Ophelia can force a draw in \ttt{} by Theorem~\ref{MB to ttt}.  
This is a design that previous techniques do not assist us with, and so in order to obtain this result, we develop some techniques for general hypergraphs in Subsection~\ref{subs:MB_HG} that significantly ease the case analysis required. 

In Section~\ref{Strong}, we consider the \ttt{} game, particularly when played on 
two types of designs, namely triple systems, and transversal designs. 
One of our most important results is Theorem~\ref{TS(v,lambda)}, where we show that playing \ttt{} on a TS$(v, \lambda)$ is a Xeno win if and only if $v \geq 5$.
In Theorem~\ref{TD(2,n)}, we show that a \TD$(2,n)$ is a Xeno win if and only if $n\geq 2$, and in Theorem~\ref{k=3 is Xeno win} we show that  a \TD$(3,n)$ is a Xeno win if and only if $n\geq 3$. We give some remarks on \TD$(4,n)$s, where we show that Ophelia can force a draw during the \ttt{} game for $n \leq 5$. We also show a surprising result in Theorem \ref{pan-hamiltonian} that answers Question 1.7 in the negative. 

In the final section we summarize our results and provide some concluding remarks.

\section{The \MB{} Game}
\label{Weak Breaker win}

In this section, we consider the \MB{} game. 
We first briefly discuss the ramifications of Theorem~\ref{thm:MB_general_breaker} and Theorem~\ref{MB win condition} on BIBDs, before moving on to general hypergraphs, which will then be used to study the game on transversal designs. 
Note that general designs are examples of $k$-uniform linear hypergraphs, and so Theorem \ref{thm:MB_general_breaker} can be applied to these structures. In Section~\ref{subs:MB_HG}, we discuss playing \MB{} on arbitrary hypergraphs. We conclude the section by applying these hypergraph methods to prove that the smallest unknown case, a \TD$(4,4)$, is Breaker win in the \MB{} game, and hence Ophelia can force a draw in the corresponding \ttt{} game.

\subsection{The \MB{} Game on BIBDs}
\label{subs:MB_BIBD}

In this section we consider the \MB{} game on BIBDs. Recall that an admissable triple $(v,k,\lambda)$ is one that satisfies the necessary conditions of Equation~(\ref{bibd nec}).

\begin{lemma}
\label{lem:BIBD MB}
    Given an admissible $(v,k,1)$, a $\mathrm{BIBD}(v,k,1)$ is Maker win if \[v > k(k-1)2^{k-3}+1,\] and is Breaker win if 
    \[
    v < \frac{-k+1 + \sqrt{(k+1)^2 + k(k-1)2^{k+2}}}{2}.
    \]
\end{lemma}
\begin{proof}
	Noting that $|\mathcal{B}|= \frac{v(v-1)}{k(k-1)}$ and each point is in $r=\frac{v-1}{k-1}$ blocks, the Maker win case follows from Theorem \ref{MB win condition}, and the Breaker win case follows from Theorem~\ref{thm:MB_general_breaker}. 
\end{proof}

\begin{corollary}\label{cor:BIBD(v,3,1} 
    Maker can win Maker-Breaker on a $\mathrm{BIBD}(v,3,1)$ ($\mathrm{STS}(v)$) if and only if $v\equiv 1, 3\bmod 6$ and $v\geq 7$. 
\end{corollary}
\begin{proof}
    By Equation~(\ref{bibd nec}), $v\equiv 1, 3\bmod 6$ is necessary for the existence of the design. 
    The design is clearly a Breaker win if $v=3$. 
    It is a Maker win if $v=7$ by Theorem \ref{Affine and Projective Planes Result - Maker Breaker}, since then it is equivalent to a $\Pi_2$.
    It is Maker win if $v > 7$ by Lemma~\ref{lem:BIBD MB}.    
\end{proof}

We now provide a complete solution in the case when $k=4$.

\begin{theorem}
\label{BIBD(v,4,1)}
Maker can win Maker-Breaker on a $\mathrm{BIBD}(v,4,1)$  if and only if $v\equiv 1, 4\bmod {12}$ and $v > 13$ (so Breaker can win on this design for $v \in \{4,13\}$).
\end{theorem}
\begin{proof}
We first note that $v\equiv 1, 4\bmod {12}$ is necessary for the existence of a BIBD$(v,4,1)$ by Equation~(\ref{bibd nec}).
    This design is Breaker win for $v\leq 12$ by Lemma~\ref{lem:BIBD MB}, but the only admissible value in this range is $v=4$. When $v=13$, the design is equivalent to $\Pi_3$ and so is a Breaker win by Theorem \ref{Affine and Projective Planes Result - Maker Breaker}. When $v=16$, the design is equivalent to $\pi_4$, so is a Maker win by Theorem \ref{Affine and Projective Planes Result - Maker Breaker}.
    When $v > 25$, the design is a Maker win by Lemma~\ref{lem:BIBD MB}. 
    The only remaining admissible value is $v=25$. 
    There are 18 non-isomorphic $\mathrm{BIBD}(25,4,1)$ (see \cite{Handbook}), and using a computational search, we have found that each of these is Maker win. 
\end{proof} 

\subsection{General Hypergraphs}
\label{subs:MB_HG}

In this section, we analyze when Breaker can win the \MB{} game played on hypergraphs in general.  
We emphasize that when Breaker can win in the \MB{} game on some hypergraph, Ophelia can force a draw in the \ttt{} game on the same hypergraph. 

Consider the \MB{} game played on arbitrary hypergraphs.
As we play the game, each Maker or Breaker move can be seen as deleting certain vertices and hyperedges from the hypergraph. 
Maker playing on a vertex corresponds to deleting that vertex from each hyperedge that contains it and deleting the vertex from the vertex set. 
Breaker playing on a vertex corresponds to deleting all hyperedges containing that vertex and removing the vertex from the vertex set. 
Thus as the game progresses, the game corresponds to a succession of nested hypergraphs until either there is an empty hyperedge in the nested hypergraph, or the nested hypergraph has no hyperedges. 
That is, either Maker has played on all vertices of a hyperedge in the original hypergraph and so Maker has won, or Breaker has played on at least one vertex of each hyperedge of the original hypergraph and so Breaker has won. 
We now state an obvious consequence of these observations. The analogous result for \ttt{} follows from Theorem~\ref{MB to ttt}.
\begin{lemma} \label{lem:del_one_move}
    Suppose we are playing \MB{} on a hypergraph $H=(V,\cE)$, and that Maker starts by playing on $u\in V$ and Breaker responds by playing on $v \in V\setminus\{u\}$. Breaker can win in the remainder of the game on $H$ if she can win on the hypergraph $H'=(V', \cE')$, where $V' = V\setminus \{u,v\}$ and $\cE' =  \{h \setminus \{u\} : h \in \cE \text{ with }  v\notin h\}$. 
\end{lemma}
Therefore, throughout the upcoming case analysis, if we have assumed that Maker played $u$ and then Breaker played $v$, we can continue the analysis on the set of hyperedges $\{h\setminus \{u\} : h \in \cE,  v\notin h\}$.
To illustrate this, suppose we wish to know if Breaker can win when playing on the hyperedges $\{a,b,c\}$, $\{a,d,e\}$, and $\{b,e,f\}$. If Maker plays $a$ and Breaker plays $f$, then we can now just ask whether Breaker can force a draw on $\{b,c\}, \{d,e\}$, having removed $a$ from each hyperedge and deleted the hyperedge $\{b,e,f\}$ since it contains $f$. 
Clearly, Breaker can win in this case. 
We can repeat this analysis for each of Maker's initial choices, where we provide Breaker's response to this move and then show that Breaker can win on the resulting hypergraph $H'$. 
As such, we may conclude that Breaker can always win on $H$. 

If we are playing either \MB{} or \ttt{} on a hypergraph $H = (V, {\cal E})$, we may use an automorphism (a map $\phi: V\to V$ which preserves hyperedges) of $H$ to help us analyze the game. 
In particular, we note that a game played with Maker moves $(X_1, X_2, \ldots)$ and with Breaker moves $(O_1, O_2, \ldots)$ will have the same outcome as the permuted strategies where Maker played $(\phi(X_1), \phi(X_2), \ldots)$ and Breaker played $(\phi(O_1), \phi(O_2), \ldots)$. 

This follows for Maker since for any hyperedge $E=\{e_1,\ldots,e_k\}\in {\cal E}$ that is completely filled by Maker (causing a win) with the original strategies, the hyperedge, $\phi(E) = \{\phi(e_1),\ldots,\phi(e_k)\}$ will be completely filled by Maker (causing a win) with the permuted strategy. 
So Maker completes a hyperedge in the original if and only if they complete a hyperedge in the permuted version.
Conversely, if Breaker plays on a hyperedge $E\in {\cal E}$, then they will have played on the hyperedge $\phi(E)$,\, and so if all hyperedges have been played on by Breaker in the original hypergraph, they will also be played in the permuted version. 
We note that a similar argument holds for the regular \ttt{} game as well.

Given two hypergraphs $H_1$ and $H_2$, define $H_1 \disunion H_2$ as the vertex-disjoint union of the two hypergraphs. 
\begin{lemma} \label{lem:disjoint_hyper}
If both $H$ and $H'$ are vertex-disjoint Breaker win hypergraphs, then Breaker can win on $H \disunion H'$. 
\end{lemma}

\begin{proof}
    On each turn while playing on $H \disunion H'$, if Maker plays on a vertex in $V(H)\subseteq V(H \disunion H')$, then Breaker plays on $V(H)$ as if she was playing her strategy on $H$, and similarly for $V(H')$. 
\end{proof}

\subsection{The \MB{} Game on Transversal Designs}
\label{subs:MB_TD}

Transversal designs are, in some sense, a relaxation of finite geometries, and so given the result of Theorems~\ref{Affine and Projective Planes Result} and \ref{Affine and Projective Planes Result - Maker Breaker}, it is natural to ask what results can be found for transversal designs.
Recall that a \TD$(k,n)$ contains $kn$ points, $n^2$ blocks, and each point is in $n$ blocks. 

 Applying Theorem~\ref{thm:MB_general_breaker} and Theorem~\ref{MB win condition} to transversal designs, we have the following corollary.  
\begin{corollary}
\label{TD 2^k > n(n+1)}
Breaker can win \MB{} on a \TD$(k,n)$ whenever $n \leq 2^{k/2}-1/2$, and Maker can win whenever $n > k\,2^{k-3}$. 
\end{corollary}
\begin{proof}
Theorem~\ref{thm:MB_general_breaker} shows that Breaker can win when $2^k > n(n+1)$. This holds when $n \leq 2^{k/2}-1/2$, which gives the first result. 
Theorem~\ref{MB win condition} shows that Maker can win when $2^{k-3}kn < n^2$, which gives the Maker win result.
\end{proof}

In particular, we have the following for small values of $k$. 
\begin{corollary}\label{cor:small values of k} 
    Maker can win \MB{} on a \TD$(3,n)$ when $n \geq 4$. 
    Maker can win \MB{} on a \TD$(4,n)$, if $n \geq 9$, and Breaker can win if $n= 3$.
\end{corollary}

We will see later in Theorem~\ref{TD(3,n)} that
Ophelia can force a draw when playing \ttt{} on a TD$(3,n)$, and thus by  Theorem~\ref{MB to ttt} 
a \TD$(3,3)$ is Maker win.
This means that the smallest unresolved case is a \TD$(4,4)$.
We thus now begin to work towards showing that playing \MB{} on a \TD$(4,4)$ is Breaker win. 
To aid us in this, we will use the hypergraph theory developed in the last subsection to show that Breaker can win on the seven hypergraphs shown in Figure~\ref{fig:graphs5}. 

\begin{lemma}
	\label{lem:drawing_hypergraphs}
	In each of the following seven hypergraphs, Breaker can win.
\begin{center}

    \begin{tabular}{c|l}
    \mbox{\rm Label}& \mbox{ \rm Hyperedges} \\ \hline
         $H_1$ & $x_1x_3x_4$, $x_3x_5u$, $x_1x_2v$, $uv$ \\
         $H_2$ & $x_1x_2x_3$, $x_3x_4u$, $x_5x_6v$, $uv$ \\
         $H_3$ & $x_1x_2x_3$, $x_3x_4x_5$, $x_5x_6 u$, $ux_7$ \\
         $H_4$ & $x_1x_2x_3$, $x_3x_4x_5$, $x_5x_6$ \\
         $H_5$ & $x_1x_2x_3$, $x_3x_4$ \\
         $H_6$ & $x_1x_2x_3$ \\
         $H_7$ & $x_1x_2$ 
    \end{tabular}
\end{center}
\end{lemma}
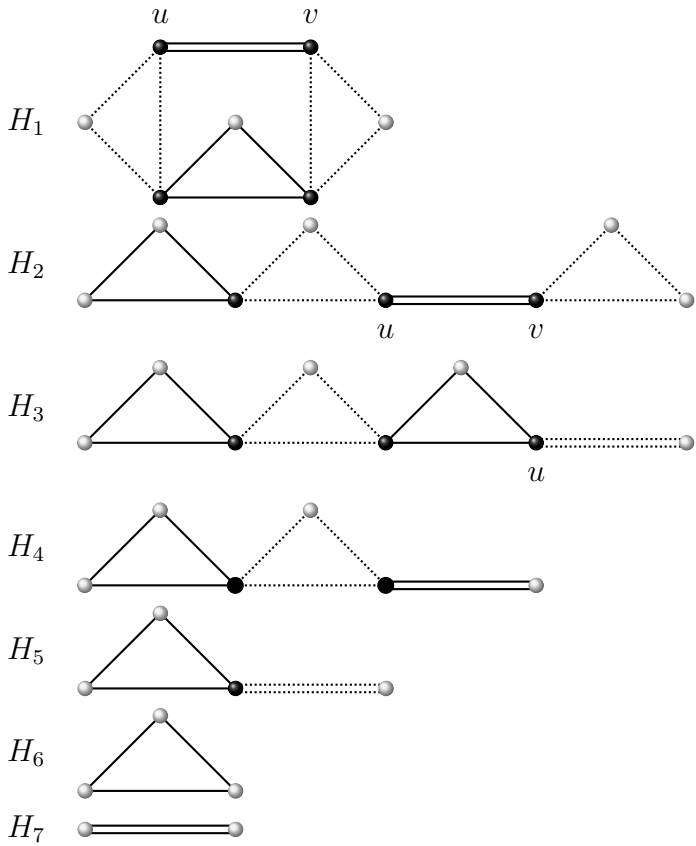
\begin{figure}
\begin{tabular}{c l c}
\raisebox{1cm}{$H_1$}& 

\begin{tikzpicture}[x=2cm,y=1cm,scale=1]
\draw[thick] (0,0) -- (1,0);
\draw[thick] (0,2.05) -- (1,2.05);
\draw[thick] (0,1.95) -- (1,1.95);
\draw[densely dotted, thick] (1,0) -- (1,2);
\draw[densely dotted, thick] (0,0) -- (0,2);
\draw[thick] (0,0) -- (0.5,1); 
\draw[thick] (0.5,1) -- (1,0); 
\draw[densely dotted, thick] (1,0) -- (1.5,1); 
\draw[densely dotted, thick] (1,2) -- (1.5,1); 
\draw[densely dotted, thick] (0,0) -- (-0.5,1); 
\draw[densely dotted, thick] (0,2) -- (-0.5,1); 
			
\shade [ball color=black!10] (-0.5,1) circle (3pt);
\shade [ball color=black!10] (1.5,1) circle (3pt);
\shade [ball color=black!10] (0.5,1) circle (3pt);

\shade [ball color=black] (0,2) circle (3pt);
\shade [ball color=black] (1,0) circle (3pt);
\shade [ball color=black] (1,2) circle (3pt);
\shade [ball color=black] (0,0) circle (3pt);

\node [above=5pt] (v) at (1,2) {$v$};
\node [above=5pt] (u) at (0,2) {$u$};

\end{tikzpicture}  \\

       \raisebox{1cm}{$H_2$}  & \begin{tikzpicture}[x=2cm,y=1cm,scale=1, 
       singlenode/.style={circle, draw=black, fill=red!50, very thick, minimum size=3pt},
       doublenode/.style={circle, draw=black, fill=red!50, very thick, minimum size=3pt}
       ]
\draw[thick] (0,0) -- (1,0);
\draw[densely dotted, thick] (1,0) -- (2,0);
\draw[thick] (2,0.05) -- (3,0.05);
\draw[thick] (2,-0.05) -- (3,-0.05);
\draw[densely dotted, thick] (3,0) -- (4,0); 
\draw[thick] (0,0) -- (0.5,1); 
\draw[thick] (0.5,1) -- (1,0); 
\draw[densely dotted, thick] (1,0) -- (1.5,1); 
\draw[densely dotted, thick] (1.5,1) -- (2,0); 
\draw[densely dotted, thick] (3,0) -- (3.5,1); 
\draw[densely dotted, thick] (3.5,1) -- (4,0); 
			
\shade [ball color=black!10] (0,0) circle (3pt);
\shade [ball color=black!10] (0.5,1) circle (3pt);
\shade [ball color=black!10] (1.5,1) circle (3pt);
\shade [ball color=black!10] (3.5,1) circle (3pt);
\shade [ball color=black!10] (4,0) circle (3pt);
\shade [ball color=black] (1,0) circle (3pt);
\shade [ball color=black] (2,0) circle (3pt);
\shade [ball color=black] (3,0) circle (3pt);
\node [below=5pt] (u) at (2,0) {$u$};
\node [below=5pt] (v) at (3,0) {$v$};
\end{tikzpicture}  \\

    \raisebox{1cm}{$H_3$} & 
    \begin{tikzpicture}[x=2cm,y=1cm,scale=1]
\draw[thick] (0,0) -- (1,0);
\draw[densely dotted, thick] (1,0) -- (2,0);
\draw[thick] (2,0) -- (3,0);
\draw[densely dotted, thick] (3,0.05) -- (4,0.05); 
\draw[densely dotted, thick] (3,-0.05) -- (4,-0.05); 
\draw[thick] (0,0) -- (0.5,1); 
\draw[thick] (0.5,1) -- (1,0); 
\draw[densely dotted, thick] (1,0) -- (1.5,1); 
\draw[densely dotted, thick] (1.5,1) -- (2,0); 
\draw[thick] (2,0) -- (2.5,1); 
\draw[thick] (2.5,1) -- (3,0); 
			
\shade [ball color=black!10] (0,0) circle (3pt);
\shade [ball color=black!10] (0.5,1) circle (3pt);
\shade [ball color=black!10] (1.5,1) circle (3pt);
\shade [ball color=black!10] (2.5,1) circle (3pt);
\shade [ball color=black!10] (4,0) circle (3pt);

\shade [ball color=black] (1,0) circle (3pt);
\shade [ball color=black] (2,0) circle (3pt);
\shade [ball color=black] (3,0) circle (3pt);
\node [below=5pt] (v) at (3,0) {$u$};
\end{tikzpicture}\\

    \raisebox{.5cm}{$H_4$} & 
    \begin{tikzpicture}[x=2cm,y=1cm,scale=1]
\draw[thick] (0,0) -- (1,0);
\draw[densely dotted, thick] (1,0) -- (2,0);
\draw[thick] (2,0.05) -- (3,0.05);
\draw[thick] (2,-0.05) -- (3,-0.05);
\draw[thick] (0,0) -- (0.5,1); 
\draw[thick] (0.5,1) -- (1,0); 
\draw[densely dotted, thick] (1,0) -- (1.5,1); 
\draw[densely dotted, thick] (1.5,1) -- (2,0); 
			
\shade [ball color=black!10] (0,0) circle (3pt);
\shade [ball color=black!10] (0.5,1) circle (3pt);
\shade [ball color=black!10] (1.5,1) circle (3pt);
\shade [ball color=black!10] (3,0) circle (3pt);

\filldraw[fill=black] (1,0) circle (3pt);
\filldraw[fill=black] (2,0) circle (3pt);
\end{tikzpicture}  \\

\raisebox{.5cm}{$H_5$} & 
    \begin{tikzpicture}[x=2cm,y=1cm,scale=1]
\draw[thick] (2,0) -- (3,0);
\draw[densely dotted, thick] (3,0.05) -- (4,0.05); 
\draw[densely dotted, thick] (3,-0.05) -- (4,-0.05); 
\draw[thick] (2,0) -- (2.5,1); 
\draw[thick] (2.5,1) -- (3,0); 
			
\shade [ball color=black!10] (2.5,1) circle (3pt);
\shade [ball color=black!10] (4,0) circle (3pt);
\shade [ball color=black!10] (2,0) circle (3pt);

\shade [ball color=black] (3,0) circle (3pt);
\end{tikzpicture}\\

\raisebox{.5cm}{$H_6$} & 
    \begin{tikzpicture}[x=2cm,y=1cm,scale=1]
\draw[thick] (2,0) -- (3,0);
\draw[thick] (2,0) -- (2.5,1); 
\draw[thick] (2.5,1) -- (3,0); 
			
\shade [ball color=black!10] (2.5,1) circle (3pt);
\shade [ball color=black!10] (2,0) circle (3pt);

\shade [ball color=black!10] (3,0) circle (3pt);
\end{tikzpicture}\\

    $H_7$ & 
    \begin{tikzpicture}[x=2cm,y=1cm,scale=1]

\draw[thick] (2,0.05) -- (3,0.05);
\draw[thick] (2,-0.05) -- (3,-0.05);

\shade [ball color=black!10] (3,0) circle (3pt);
\shade [ball color=black!10] (2,0) circle (3pt);
\end{tikzpicture}  

\end{tabular}
\caption{Seven hypergraphs, with a hyperedge of cardinality $2$ indicated as a double line and of cardinality $3$ indicated as a triangle (either solid or dashed). Each vertex that is incident with one hyperedge is white, and with two hyperedges is black.}
    \label{fig:graphs5}
\end{figure}
\begin{proof}
	These hypergraphs are displayed in Figure \ref{fig:graphs5}. 
	Recall that by Theorem~\ref{score<1}, if the score of the game is below $1$ after Maker takes his turn, then Maker cannot win.     
	The score of each of these hypergraphs is at most $5/8$.  
	As such, Maker must play on a vertex that increases the score by at least $3/8$.
    Each such vertex has been indicated by being labelled as $u$ or $v$. 
    Note that there are no such vertices in $H_4$, $H_5$, $H_6$, and $H_7$, and so the proof is complete in these cases.
    
    In $H_1$, if Maker chooses $u$, then Breaker responds with $v$. Lemma \ref{lem:del_one_move} completes the proof, when we note that the hypergraph 
    \[
    	H'=H'(V(H_1)\setminus \{u,v\}, \{h\setminus \{u\} : h \in E(H_1),  v\notin h\}) = H_5,
    \] 
    on which we already know that Breaker can win.
    The case that Maker chooses $v$ and Breaker responds with $u$ follows in the same way.

    In $H_2$, if Maker chooses $u$, then Breaker responds with $v$. 
    Lemma \ref{lem:del_one_move} completes the proof, noting that $H'=H_5$.
    If Maker chooses $v$, then Breaker responds with $u$. 
    Lemmas \ref{lem:del_one_move} and \ref{lem:disjoint_hyper} complete the proof, noting that $H'=H_6\disunion H_7$.

    In $H_3$, if Maker chooses $u$, then Breaker responds by playing $x_7$, the other vertex in the hyperedge of cardinality $2$. 
    Lemma \ref{lem:del_one_move} completes the proof, noting that $H'=H_4$.
\end{proof}

\begin{theorem}
	\label{TD(4,4)}  
  Breaker can win \MB{} on a \TD$(4,4)$. 
\end{theorem}
\begin{proof}
	We note that up to isomorphism, there is only one \TD$(4,4)$, see \cite{WanlessEgan_EnumMOLS}.
	The \TD$(4,4)$ has $16$ points, $16$ blocks and $4$ groups.
	Using the notation defined in \cite{DHMM}, we take the point set to be $X = \{r_i,c_i,\alpha_i,\beta_i : i \in \{1,2,3,4\}\}$. 
	Define the groups $G_1=\{r_i : i\in \{1,2,3,4\}\}$, $G_2=\{c_i : i\in \{1,2,3,4\}\}$, $G_3=\{\alpha_i : i\in \{1,2,3,4\}\}$, and $G_4=\{\beta_i : i\in \{1,2,3,4\}\}$, and note that the groups partition the point set $X=G_1\cup G_2\cup G_3\cup G_4$. 
	The set of blocks $\mathcal{B}$ are: 
	\[
	\begin{array}{lcccc}
		\{r_1, c_1, \entryi_1, \entryii_1\}, & \{r_4, c_2, \entryi_3, \entryii_4\}, & \{r_2, c_3, \entryi_4, \entryii_2\}, & \{r_3, c_4, \entryi_2, \entryii_3\}, \\
		
		\{r_3, c_3, \entryi_1, \entryii_4\}, & \{r_1, c_2, \entryi_2, \entryii_2\}, & 
		\{r_4, c_1, \entryi_4, \entryii_3\}, & \{r_2, c_4, \entryi_3, \entryii_1\}, \\
		
		\{r_2, c_1, \entryi_2, \entryii_4\}, & \{r_3, c_2, \entryi_4, \entryii_1\}, & 
		\{r_1, c_3, \entryi_3, \entryii_3\}, &\{r_4, c_4, \entryi_1, \entryii_2\}, \\
		
		\{r_3, c_1, \entryi_3, \entryii_2\}, & \{r_2, c_2, \entryi_1, \entryii_3\}, & \{r_4, c_3, \entryi_2, \entryii_1\}, & \{r_1, c_4, \entryi_4, \entryii_4\}. \\
	\end{array}
	\]

	From now on, we will use hypergraph notation to align with the previous three lemmas. As such, we define the hypergraph at play as $H=(V,\cE)$ with $V=X$ and $\cE=\cB$. 
	As described in Section~\ref{subs:MB_HG}, when Maker plays on a point, we remove that point from $V$ and from every hyperedge containing it, and when Breaker plays on a point, we remove it and all hyperedges it contains from $H$.
    We will use the following automorphisms of $H$: 
	\begin{align*}
		\phi_1 &= (c_1 c_3)(c_2 c_4)(\alpha_1 \alpha_3)(\alpha_2 \alpha_4) (\beta_1 \beta_3)(\beta_2 \beta_4), \\
		\phi_2 &= (c_1 c_2)(c_3 c_4)(\alpha_1 \alpha_3)(\alpha_2 \alpha_4) (\beta_1 \beta_2)(\beta_3 \beta_4), \\
		\phi_3 &= (c_1 \alpha_1 \beta_1)(c_2 \alpha_3 \beta_4)(c_3 \alpha_4 \beta_2) (c_4 \alpha_2 \beta_3), \\
		\phi_4 &= (r_2 r_3) (c_3 c_4)(\alpha_1 \beta_1)(\alpha_2 \beta_2)(\alpha_3 \beta_4)(\alpha_4 \beta_3), \\
		\phi_5 &= (r_3 r_4) (c_1 \alpha_2)(c_2 \alpha_1)(c_3 \alpha_3)(c_4 \alpha_4)(\beta_1 \beta_2).
	\end{align*}
    Recall that showing a sequence of moves leads to a win for either player shows that the sequence of moves permuted by an automorphism of $H$ leads to a win for the same player.
	Given any pair of points $u$ and $v$ not in $G_1$, there exists an automorphism formed by combining $\phi_1$, $\phi_2$, and $\phi_3$ in some way, such that $u$ is mapped to $v$ but where $r_i$ is a fixed point for each $i$.
	Let $\mathcal{A}$ denote the set of these automorphisms. 
 
	We will give a strategy for Breaker for every possible Maker move.	
	To start, we note that the automorphism group of $H$ acts transitively on its vertices.  
	Therefore Maker's first move is arbitrary, and can be selected to be $r_1$. 
	Breaker's strategy will be to play her first two moves in $G_1$. 
	As such, over Maker's first three moves, he will have either one or two of those moves played on $G_1$. In the following, we break into these two cases.
	Note that we will make extensive use of Theorem~\ref{score<1}: Breaker wins if the score falls below $1$ at the end of Maker's turn.

	%-----------------------------------------
	{\bf Case 1:} Two of Maker's first three moves are in $G_1$. 
	By applying the automorphisms $\phi_4$ and $\phi_5$ in some order, we may assume that Breaker played her first two moves on $r_2$ and $r_4$ and that Maker played two of his first three moves on $r_1$ and $r_3$. Using the automorphisms in $\mathcal{A}$, we may map any vertex to $c_1$ while retaining each $r_i$ as a fixed point. We can therefore assume the remaining Maker move of Maker's first three moves was played on $c_1$. 
	
	Breaker responds by playing $\alpha_1$ for her third move. 
	There are six remaining hyperedges, namely 
	$\{\alpha_3, \beta_2\}$, $\{c_2, \alpha_2, \beta_2\}$, $\{c_2, \alpha_4, \beta_1\}$,$\{c_3, \alpha_3, \beta_3\}$,$\{c_4, \alpha_2, \beta_3\}$, and $\{c_4, \alpha_4, \beta_4\}$; see Figure \ref{fig:TD44_hyper}. 
	\begin{figure}[bth]
		\begin{center}
				\begin{tikzpicture}[x=2cm,y=0.5cm,scale=1]
				\draw[thick] (-2,2) -- (-1,2);%beta1---alpha4
				\draw[densely dotted, thick] (-1,2) -- (0,2);%alpha4---c4
				
				\draw[densely dotted, thick]  (-0.5,4) -- (0,2);%beta4---alpha4
				 \draw[densely dotted, thick]  (-0.5,4) -- (-1,2);%beta4---alpha4 (error -> c4?)
				 
				 \draw[densely dotted, thick]  (-1.5,0) -- (-0.5,0);%c2---alpha2
				\draw[densely dotted, thick]  (-0.5,0) -- (-1,-2);%alpha2---beta2
				\draw[densely dotted, thick]  (-1.5,0) -- (-1,-2);%c2---beta2
				
				 \draw[thick]  (0,2) -- (-0.5,0);%c4---alpha2
				\draw[thick]  (-0.5,0) -- (0.5,0);%alpha2---beta3
				\draw[thick]  (0,2) -- (0.5,0);%c4---beta3
				
				 \draw[thick]  (-2,2) -- (-1.5,0);%beta1---c2
				\draw[thick]  (-1.5,0) -- (-1,2);%c2---alpha2
				
				\draw[densely dotted, thick]  (0.5,0) -- (1.5,0);%beta3---c3
				\draw[densely dotted, thick]  (1.5,0) -- (1,-2);%c3---alpha3
				\draw[densely dotted, thick]  (0.5,0) -- (1,-2);%beta3---alpha3
	
	            \draw[thick] (-1,-1.95) -- (1,-1.95);
	            \draw[thick] (-1,-2.05) -- (1,-2.05);
				\node[above=4pt] (beta4) at (-0.5,4) {$\beta_4$};
				\shade [ball color=black] (-0.5,4) circle (3pt);
				
				\node [left=5pt] (beta1) at (-2,2) {$\beta_1$};
				\shade [ball color=white] (-2,2) circle (3pt);
				
				\node[=5pt](alpha4) at (-1.2,2.4) {$\alpha_4$};
				\shade[ball color=black] (-1,2) circle (3pt);
				
				\node[right=4pt] (c4) at (0,2) {$c_4$};
				\shade [ball color=black] (0,2) circle (3pt);
				
				\node[left=4pt] (c2) at (-1.5,0) {$c_2$};
				\shade [ball color=black] (-1.5,0) circle (3pt);
				
				\node[below=4pt] (alpha2) at (-0.3,0.1) {$\alpha_2$};
				\shade [ball color=black] (-0.5,0) circle (3pt);
				
				\node[below=4pt] (beta3) at (0.4,0.1) {$\beta_3$};
				\shade [ball color=black] (0.5,0) circle (3pt);
				
				\node[right=4pt] (c3) at (1.5,0) {$c_3$};
				\shade [ball color=white] (1.5,0) circle (3pt);
				
				\node[left=4pt] (beta2) at (-1,-2) {$\beta_2$};
				\shade [ball color=black] (-1,-2) circle (3pt);
				
				\node[right=4pt] (alpha3) at (1,-2) {$\alpha_3$};
				\shade [ball color=black] (1,-2) circle (3pt);
				
				\node[above=4pt] (beta4) at (-0.5,4) {$\beta_4$};
				\shade [ball color=white] (-0.5,4) circle (3pt);
			\end{tikzpicture}
			\caption{Remaining vertices of the six hyperedges in which Breaker has not yet played in Case 1 of Theorem~\ref{TD(4,4)}. These hyperedges are  $\{\alpha_3, \beta_2\}$, $\{c_2, \alpha_2, \beta_2\}$, $\{c_2, \alpha_4, \beta_1\}$,$\{c_3, \alpha_3, \beta_3\}$,$\{c_4, \alpha_2, \beta_3\}$, and $\{c_4, \alpha_4, \beta_4\}$.}
			\label{fig:TD44_hyper}
		\end{center}
	\end{figure}
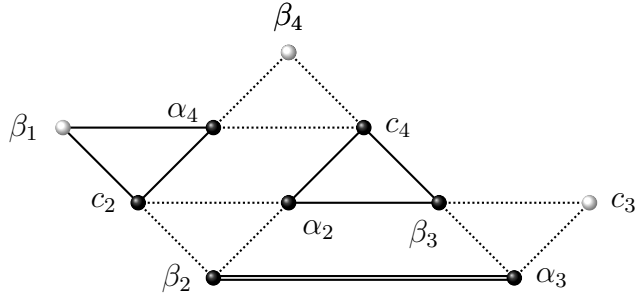
	It is now Maker's move, and he has ten possible choices for his next move.
    Table \ref{tab:last_moves} gives the remainder of the game. For each of Maker's possible moves, we give Breaker's response and the hypergraph that results. In each of these resulting hypergraphs, Breaker can win by Lemmas~\ref{lem:disjoint_hyper} and \ref{lem:drawing_hypergraphs}.
	%-----------------------------------------
	\begin{table}[thb] 
		\begin{center}
			\renewcommand{\arraystretch}{1.25}
			\begin{tabular}{c c}
			$\begin{array}{|c|c|c|c|c|} \hline
			X & O & \text{Remaining hypergraphs} \\ \hline
			c_2 &  \beta_2 & H_3 \\ \hline
			c_3 &  \beta_3 & H_3 \\ \hline
			c_4 &  \alpha_2 & H_5 \disunion H_5 \\ \hline
			\alpha_2 &  \beta_2 & H_2 \\ \hline
			\alpha_3 &  \beta_2 & H_3 \\ \hline
			\end{array}$
			$\begin{array}{|c|c|c|c|c|} \hline
			X & O& \text{Remaining hypergraphs} \\ \hline
			\alpha_4 &  \beta_2 & H_4 \disunion H_7 \\ \hline
			\beta_1 &  \beta_2 & H_3 \\ \hline
			\beta_2 &  \alpha_3 & H_1 \\ \hline
			\beta_3 &  \alpha_3 & H_1 \\ \hline
			\beta_4 &  \beta_2 & H_2 \\ \hline
			\end{array}$
			\end{tabular}
			\caption{The fourth Maker move after each player has played three moves in Case $1$, Breaker's response, and the remaining hypergraph.}
			\label{tab:last_moves}
		\end{center}
	\end{table}

%-----------------------------------------
	{\bf Case 2:} Exactly one of Maker's first three moves is in $G_1$.  
	By applying the automorphism $\phi_4$ and $\phi_5$ in some order, we can assume Breaker played her first two moves on $r_2$ and $r_3$. By using an automorphism in $\mathcal{A}$,  we can assume that Maker plays his second move on $c_1$.  We break into subcases. 
 %-----------------------------------------
 
	{\bf Subcase 2a:} Maker's third move is in the remaining hyperedge $\{\alpha_1, \beta_1\}$.
	
	Note that automorphism $\phi_4$ transposes $\beta_1$ and $\alpha_1$, the two possibilities for Maker's third move; transposes $r_2$ with $r_3$, which both contain Breaker moves, and fixes both $r_1$ and $c_1$, which contain Maker moves. We can therefore assume without loss of generality that Maker played on $\alpha_1$. 
	In response, Breaker plays her third move on $\beta_1$, to stop Maker from winning on the next turn. 
	The hyperedges that Breaker has not played on are the five hyperedges that contain three vertices unoccupied by Maker,  
	$\{c_2,\alpha_2, \beta_2\}$, 
	$\{c_3,\alpha_3, \beta_3\}$, 
	$\{c_4,\alpha_4, \beta_4\}$, 
	$\{r_4,\alpha_4, \beta_3\}$, 
	and $\{r_4, c_4, \beta_2\}$; and one hyperedge with no vertices occupied by Maker, $\{r_4, c_2, \alpha_3, \beta_4\}$. 
	As a result, the current score is $\frac{11}{16}$. The only way that Maker can increase this score to at least $1$ on his turn is to play on $r_4$, so he does this. 
	In response, Breaker plays $\beta_3$, reducing the score to $\frac{5}{8}$. The remaining hyperedges are $\{c_2,\alpha_2, \beta_2\}$, $\{c_4,\alpha_4, \beta_4\}$, $\{c_4, \beta_2\}$, and $\{c_2, \alpha_3, \beta_4\}$. But this is isomorphic to the hypergraph $H_1$ from Lemma~\ref{lem:drawing_hypergraphs}, and so Breaker can win.

	%-----------------------------------------
	{\bf Subcase 2b:} Maker's third move is not in $\{\alpha_1, \beta_1\}$.
	
	Recall that Maker's moves so far are $\{r_1,c_1\}$ and Breaker's are $\{r_2,r_3\}$. 
	Whatever Maker's third move in this subcase, no hyperedge contains three Maker moves since Maker's third move is not in $\{\alpha_1, \beta_1\}$, and so Breaker's turn is not forced. 
	Since Maker did not play his second or third turn in $G_1$, $r_4$ has not been played at this point. 
    Breaker plays her third move as $r_4$. 
    Now there are exactly four hyperedges that have not been played on by Breaker, which are those hyperedges that originally contained the vertex $r_1$. 
    These hyperedges are pairwise disjoint as $r_1$ was the only common vertex of these hyperedges before $r_1$ was removed by Maker playing on it, and so for the rest of play, if Maker plays on one of these hyperedges, there is always at least one remaining unplayed vertex on this hyperedge. 
    Then Breaker will respond by playing on that unplayed vertex.
    This continues for at most four turns, at which point Breaker has played on at least one vertex in each hyperedge, and so has won. 
\end{proof}

\section{The \ttt{} Game}
\label{Strong}

In this section, we consider the \ttt{} game, particularly when played on designs. We first summarize the ramifications of the results from Section~\ref{Weak Breaker win} in light of Theorem~\ref{MB to ttt}. 
We then give a complete solution for triple systems in Section~\ref{Triple Systems}. Finally in Section~\ref{TD} we give results for \TD$(k,n)$ with $k\in \{2,3\}$.

Theorem~\ref{MB to ttt} says that if Breaker can win \MB{}, then Ophelia can force a draw in \ttt{}. 
Hence we have a number of results for \ttt{} that follow from those for \MB{} from Section~\ref{Weak Breaker win}. We summarize these in the theorem below.
\begin{theorem}
\label{thm:MB to ttt}
    Ophelia can force a draw when playing \ttt{} on the following designs:
    \begin{itemize}
        \item 
        \BIBD$(v, k, 1)$ whenever $v < \frac{-k+1 + \sqrt{(k+1)^2 + k(k-1)2^{k+2}}}{2}$; \hfill {\rm (Lemma~\ref{lem:BIBD MB})}
        \item  \TD$(k,n)$, whenever $n \leq 2^{k/2}-1/2$; and \hfill {\rm (Corollary~\ref{TD 2^k > n(n+1)})}
        \item \TD$(4,4) $\hfill {\rm (Theorem~\ref{TD(4,4)})}
    \end{itemize}
\end{theorem}

\subsection{Triple Systems}
\label{Triple Systems}

Suppose that there are two non-equal blocks in a TS$(v,\lambda)$ intersecting in a point $x$, $B_1=\{x,a,b\}$ and $B_2=\{x,c,d\}$ with $b\neq d$ and that one player has played on $a$ and $c$ (note that it is possible that $a=c$). 
Further, suppose that we are in a state where neither player is about to win.
If the player who holds $a$ and $b$ now chooses $x$, they will win the game on their next move, since their opponent cannot block the wins on both $B_1$ and $B_2$ simultaneously.
We call such a situation a \emph{scissor} on $x$. 

We first deal with the case when $\lambda = 1$ and give the following theorem for playing \ttt{} on a STS$(v) = $ TS$(v,1)$, which was proved independently in~\cite{Abbey}. 

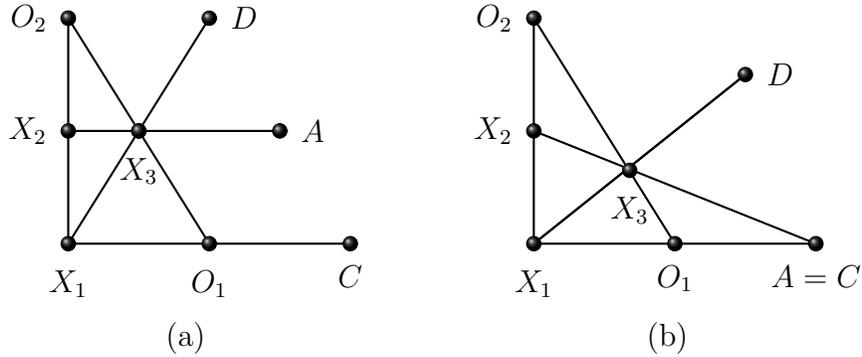
\begin{figure}[th]
	\begin{center}
		\begin{tabular}{ccc}
			\begin{tikzpicture}[x=1.25cm,y=1cm,scale=1.5]
			\draw[thick] (0,0) -- (2,0);
			\draw[thick] (0,0) -- (0,2);
			\draw[thick] (0,2) -- (1,0);
			\draw[thick] (0,0) -- (1,2);
			\draw[thick] (0,1) -- (1.5,1);
			
			\node [below=6pt] (X1) at (0,0) {$X_1$};
			\shade [ball color=black] (0,0) circle (2pt);
			
			\node[below=6pt] (O_1) at (1,0) {$O_1$};
			\shade[ball color=black] (1,0) circle (2pt);
			\shade[ball color=black] (2,0) circle (2pt);
			
			\node[left=4pt] (X2) at (0,1) {$X_2$};
			\shade [ball color=black] (0,1) circle (2pt);
			\shade [ball color=black] (0,2) circle (2pt);
			
			\node[left=4pt] (O2) at (0,2) {$O_2$};
			\shade [ball color=black] (0,2) circle (2pt);
			\shade [ball color=black] (0.5,1) circle (2pt);
			
			\node[below=6pt] (X3) at (.5,1) {$X_3$};
			\shade  [ball color=black] (.5,1) circle (2pt);
			\shade [ball color=black] (1,2) circle (2pt);
			\shade [ball color=black] (1.5,1) circle (2pt);
			\node [right=4pt] (D) at (1,2) {$D$};
			\node [right=4pt] (A) at (1.5,1) {$A$};
			\node [below=4pt] (C) at (2,0) {$C$};
			\end{tikzpicture}
			& \hspace{2ex} &
			\begin{tikzpicture}[x=1.25cm,y=1cm,scale=1.5]
			\node (X_1) at (0,0) {};
			\node (X_2) at (0,1) {};
			\node (O_2) at (0,2) {};
			\node (O_1) at (1,0) {};
			\node (X_3) at (.68,.65) {};
			\node (D) at (1.5,1.5) {};
			\node (A) at (2,0) {};
			
			\draw[thick] (0,0) -- (0,2); %(X_1) -- (O_2);
			\draw[thick] (0,0) -- (2,0); %(X_1) -- (A);
			\draw[thick] (0,0) -- (1.5,1.5); %(X_1) -- (D);
			\draw[thick] (1,0) -- (0,2); %(O_1) -- (O_2);
			\draw[thick] (0,0) -- (1.5,1.5); %(X_1) -- (D);
			\draw[thick] (0,1) -- (2,0); %(X_2) -- (A);
			
			\node [below=6pt] at (X_1) {$X_1$};
			\shade [ball color=black] (X_1) circle (2pt);
			
			\node[below=4pt]  at (O_1) {$O_1$};
			\shade[ball color=black] (O_1) circle (2pt);
			
			\node[left=4pt] at (X_2)  {$X_2$};
			\shade [ball color=black] (X_2) circle (2pt);
			
			\node[left=4pt] at (O_2) {$O_2$};
			\shade [ball color=black] (O_2) circle (2pt);
			
			\node[below=6pt] at (X_3) {$X_3$};
			\shade  [ball color=black] (X_3) circle (2pt);
			
			\node [right=4pt] at (D)  {$D$};
			\shade  [ball color=black] (D) circle (2pt);
			
			\node [below=4pt] at (A) {$A=C$};
			\shade  [ball color=black] (A) circle (2pt);
			\end{tikzpicture}
			\\
			(a) & & (b)
		\end{tabular}		
	\caption{
		Scissor on $X_3$, after the first three Xeno moves played on an STS$(v)$. 
		It is possible that $A=C$, in which case we have the configuration shown on the right.}
	\label{fig:STS($n$)}
	\end{center}
\end{figure}

\begin{theorem}
\label{STS(v)}
	Xeno can win \ttt{} on an {\rm STS}$(v)$ if and only if $v\equiv 1,3\bmod 6$ and $v > 3$.
\end{theorem}
\begin{proof}
	We first note that $v\equiv 1,3\bmod 6$ is necessary for the existence of an STS$(v)$ by Equation~(\ref{bibd nec}) from Subsection~\ref{subs Designs}.  
	It is clear that Ophelia can force a draw when $v=3$, as she plays on the only block of the design on her first turn, preventing Xeno from winning.

	Now assume that $v\geq 7$, as this is the next value when a $\mathrm{STS}(v)$ exists. 
	We provide a winning strategy for Xeno. 
	Xeno chooses his first move, $X_{1}$, arbitrarily. 
	Now, wherever Ophelia takes her first move $O_1$, there is a block $B = \{X_{1}, O_{1},C\}$ of the design; see Figure~\ref{fig:STS($n$)}. Xeno now chooses his second move $X_{2},$ where $X_{2}$ is not in $B$.
	After Xeno's second  move, Ophelia is threatened  in the block which contains Xeno's first two moves.
	So, Ophelia is forced to choose her second move, $O_{2}$, in this block.  
	Now, Xeno is threatened in the block containing Ophelia's first two moves, so he must choose his third move, $X_{3},$ in this block. 
	Now Xeno has a scissor on $X_3$, with the blocks $\{X_3, X_1,D\}$ and $\{X_3, X_2,A\}$, and so wins the game on his next move.
\end{proof}

Note that depending on the structure of the design and the points chosen, it is possible that $A=C$, so the blocks form a Pasch configuration, see Figure~\ref{fig:STS($n$)}(b), but this does not affect the argument.

We now completely solve the problem of playing \ttt{} on a Triple System, TS($v,\lambda$). We begin with the following lemma, which deals with cases where $v$ is small. 

\begin{lemma}
	\label{TS v<5}
	Ophelia can force a draw when playing \ttt{} on {\rm TS}$(v, \lambda)$ with $v < 5$, but TS$(5,\lambda)$ is a Xeno win.
\end{lemma}
\begin{proof}
	The result is immediate when $v=3, 4$. In both cases each player has fewer than three moves and so neither can complete a block, leading to a draw.
	
	We note that a TS$(5,\lambda)$ must have $\lambda$ divisible by 3 and consists of all possible triples of $K_5$ repeated $\lambda/3$ times. Thus if $K_5$ is labelled with $\{0,1,2,3,4\}$ then, without loss of generality, play starts as $X_1 = 0$, $O_1=1$, and $X_2 = 2$.
    Xeno now has a scissor on $X_2=2$ and can win by either playing 3 or 4, on the blocks $\{0,2,3\}$ or $\{0,2,4\}$, and Ophelia cannot block both.
\end{proof}

Given a TS$(v,\lambda)$, $(X, \cB)$, and a point $x\in X$, we define $G[x]$ to be the multigraph (parallel edges allowed) on the vertex set $X \setminus \{x\}$ and containing an edge $yz$ for each block $\{x,y,z\}\in \cB$. 
We note that $G[x]$ is $\lambda$-regular and saying that the edge $yz$ appears with multiplicity $\mu$ is equivalent to saying that the block $\{x,y,z\}$ appears $\mu$ times in $\cB$.	

Let $C$ be a component of $G[x]$ of order greater than 2, and suppose that there is a vertex $z \in C$ with two neighbours $w_1,w_2 \in C$. 
If at any point, Xeno has played on both $x$ and $z$, and Ophelia has not played on $w_1$ or $w_2$, then Xeno has a scissor on $z$, and so can win on the next turn by playing either $w_1$ or $w_2$. 

\begin{example}
 \label{ex:ts62}
Playing \ttt{} on a $\text{TS}(6,2)$ is Xeno win. 
	\begin{figure}
		\begin{center}
		\begin{tikzpicture}[x=1cm,y=1cm,scale=.5]
			\draw (72:3) -- (144:3) -- (216:3) -- (288:3) -- (0:3) -- cycle;
			\foreach \x in {1,...,5} {
				\node (p\x) at (\x*360/5:3.7) {};
				\shade [ball color=black]  (\x*360/5:3) circle (6pt);
			}
			\node () at (p1) {1};
			\node () at (p2) {3};
			\node () at (p3) {5};
			\node () at (p4) {4};
			\node () at (p5) {2};
		\end{tikzpicture}
			\caption{\label{fig:G_0} The (multi)graph $G[0]$ derived from the TS$(6,2)$ of Example \ref{ex:ts62}, which is the $C_5$ that results from including an edge $xy$ for each block $\{0,x,y\}$ in the TS.}
		\end{center}
	\end{figure}
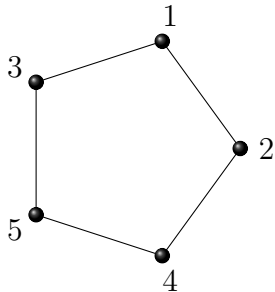
    
	 Up to isomorphism \cite{Handbook}, the blocks of a TS$(6,2)$ are
	\[
	\begin{array}{l}
	\{0,1,2\}, \{0,1,3\}, \{0,2,4\}, \{0,3,5\}, \{0,4,5\},\\
	 \{1,2,4\}, \{1,4,5\}, \{2,3,4\}, \{2,3,5\} .
	\end{array}
	\]				
    We play \ttt{} on such a TS$(6,2)$ and show that it is a Xeno win. 
    
    Xeno plays $X_1 = 0$ as his first move and note that $G[0]$ is a cycle of length $5$, $(12453)$, see Figure~\ref{fig:G_0}. If Ophelia plays any point of the design that corresponds to a vertex on this cycle, Xeno responds by playing on a point corresponding to the vertex that is two steps around the cycle, creating a scissor. For example, if Ophelia plays $O_1 = 1$, Xeno plays $X_2 = 4$ and now has created a scissor on $4$, and so can win on either of the blocks $\{0, 4, 2\}$ or $\{0, 4, 5\}$ on his next turn, since Ophelia cannot block both. 
\end{example}

\begin{theorem}	
	\label{TS(v,lambda)}
    Xeno can win \ttt{} on a TS$(v, \lambda)$ if and only if $v \geq 5$.
\end{theorem}
\begin{proof}
	Lemma \ref{TS v<5} deals with the cases when $v \leq 5$,
    so we may assume that $v>5$.
	We suppose for the sake of a contradiction that for some $v>5$, $(X, \cB)$ is a TS($v, \lambda$) in which Ophelia can force a draw.
	
	If for every $x\in X$ each connected component of $G[x]$ consists of exactly two vertices, then $\cB$ consists of the blocks of an STS$(v)$ repeated $\lambda$ times, and so Xeno can win on this design by Theorem~\ref{STS(v)}. We may thus assume that there is an $x\in X$ such that $G[x]$ has a component of order at least 3. 	

	Further, we note that given a component of order at least 3 in $G[x]$, if some vertex $y$ in that component has only one neighbour in $G[x]$, say $z$, then the edge $yz$ must occur with multiplicity $\lambda$, and so $z$ also has only one neighbour in $G[x]$, contradicting the assumption that this component has order at least 3. 
	As a result, we may assume each vertex in a component of $G[x]$ of order $3$ or more has at least two neighbours in that component. 
	
	Let $C$ be a component of $G[x]$ of order $3$ or more. Assume that Xeno has played on $x$ in the design $(X,\cB)$.  
	If $C$ contains two vertices of degree $3$ or more, say $y$ and $z$, then Xeno will have a scissor on $z$ during his next turn unless Ophelia plays her turn on $z$, in which case Xeno will have a scissor on $y$. 
	Similarly, if $C$ contains a cycle of length $4$, say $(y, a_1, b, a_2)$, Xeno will have a scissor on $y$ in $(X,\cB)$ unless Ophelia plays on $a_1$, $y$, or $a_2$. But then Xeno will instead have a scissor on $a_2$, $b$, or $a_1$, respectively. 

  	\begin{figure}[ht]
		\begin{center}
			\begin{tikzpicture}[x=1cm,y=1cm,scale=0.7]
				\shade  [ball color=black] (0,2) circle (4pt);
				\foreach \x in {-2, 2} {
					\shade  [ball color=black] (\x,1) circle (4pt);
				}
				\foreach \x in {-3, 3} {
					\shade  [ball color=black] (\x,-1) circle (4pt);
				}
				\node (a1) at (-2.5, 1.5) {$a_1$};
				\node (a2) at (2.5, 1.5) {$a_2$};
				\node (y) at (0,2.5) {$y$};
				\draw (-3,-1) -- (-2,1) -- (0,2) -- (2,1) -- (3,-1);
				\draw [style=dashed, color=black] (-3,-1) -- (0,2);
				\node (z1) at (-3.4,-1) {$b_1$};
				\draw [style=dashed, color=black] (3,-1) -- (0,2);
				\node (z2) at (3.5,-1) {$b_2$};
			\end{tikzpicture}
			\caption{\label{fig:P}The path $P$ in the component $C$.}
		\end{center}
	\end{figure}
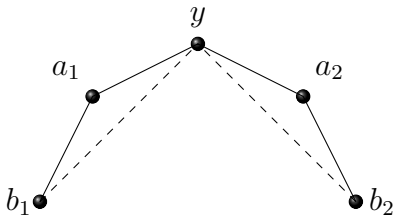

	Now suppose that $C$ contains a path $P$ of length $4$ with centre vertex $y$ and endpoints $b_1$ and $b_2$; see Figure~\ref{fig:P}.
	If Ophelia does not play on $y$ in $(X,\cB)$, Xeno can find a scissor on $a_1$ unless Ophelia played $a_1$ or $b_1$, in which case Xeno has a scissor on  $a_2$. So, Ophelia must play on $y$ during this turn.
	Now, since every vertex in $C$ connects to at least two other vertices in $C$, the endpoint $b_i$	must connect to some vertex other than $a_i$, $i=1,2$.
	If $b_i$ connects to some vertex other than $y$, then Xeno has a scissor on $b_i$ immediately before Xeno's turn, so we may assume that each $b_i$ is connected to $y$.
	Thus $C$ must look like a collection of triangles sharing a common vertex, $y$.

	Now, suppose that $y$, $z_1$, $z_2$ is a triangle in $C$, where $y$ may be connected to other vertices in $C$, but $z_1$ and $z_2$ are not, and the edges $z_1z_2$, $yz_1$ and $yz_2$ appear with multiplicity $\lambda_1$, $\lambda_2$ and $\lambda_3$ respectively in $G[x]$. Since every vertex has degree $\lambda$ in $C$, we have that 
	$$\lambda_1+\lambda_2=\lambda_1+\lambda_3 = 
	\lambda, \quad\mbox{ so }\quad \lambda_2 = \lambda_3 = \lambda - \lambda_1.$$
	Further, considering the degree of $y$ in $G[x]$, we have that $\lambda_2+\lambda_3 = 2(\lambda-\lambda_1) \leq \lambda$, and so $\lambda\leq 2\lambda_1$.

	Let $B_1=\{x,z_1,z_2\}$, $B_2=\{x,y,z_1\}$, $B_3=\{x,y,z_2\}$, and $B_4=\{y,z_1,z_2\}$. 
	The analysis above implies that the blocks $B_2$ and $B_3$ appear exactly $\lambda - \lambda_1$ times in $\cB$, and the block $B_1$ appears exactly $\lambda_1$ times in $\cB$.
	So in $G[y]$, the edges $xz_1$ and $xz_2$ must appear $\lambda-\lambda_1$ times, and since at most one vertex has more than two distinct neighbours, $z_1z_2$ has multiplicity $\lambda_1$ in $G[y]$, or else $z_1$ and $z_2$ would both have more than two distinct neighbours. Thus, the block $B_4$ also appears exactly $\lambda_1$ times in $\cB$; see Figure~\ref{fig:Multiplicity}.

	\begin{figure}
	\begin{center}
		\begin{tabular}{ccc}
			\begin{tikzpicture}[x=1cm,y=1cm,scale=.7]
			\shade [ball color=black]  (0,2) circle (4pt);
			\shade [ball color=black]  (-1,0) circle (4pt);
			\shade [ball color=black]  (1,0) circle (4pt);
			\draw (-1,0) -- (1,0) -- (0,2) -- cycle;
			\node[above] (y) at (0,2) {$y$};
			\node[left] (z1) at (-1,0) {$z_1$};
			\node[right] (z2) at (1,0) {$z_2$};
			\node[below] (z1z2) at (0,0) {$\lambda_1$};
			\node (yz1) at (-1.6,1) {$\lambda-\lambda_1$};
			\node (yz2) at (1.6,1) {$\lambda-\lambda_1$};
			\node (Gx) at (0,-1.7) {$G[x]$};
			\end{tikzpicture}
			\;\;
			&				
			\begin{tikzpicture}[x=1cm,y=1cm,scale=.7]
			\shade [ball color=black]  (0,2) circle (4pt);
			\shade [ball color=black]  (-1,0) circle (4pt);
			\shade [ball color=black]  (1,0) circle (4pt);
			\draw (-1,0) -- (1,0) -- (0,2) -- cycle;
			\node[above] (x) at (0,2) {$x$};
			\node[left] (z1) at (-1,0) {$z_1$};
			\node[right] (z2) at (1,0) {$z_2$};
			\node[below] (z1z2) at (0,0) {$\lambda_1$};
			\node (yz1) at (-1.6,1) {$\lambda-\lambda_1$};
			\node (yz2) at (1.6,1) {$\lambda-\lambda_1$};
			\node (Gy) at (0,-1.7) {$G[y]$};
			\end{tikzpicture}
			\;\;
			&
			\begin{tikzpicture}[x=1cm,y=1cm,scale=.7]
			\shade [ball color=black]  (0,2) circle (4pt);
			\shade [ball color=black]  (-1,0) circle (4pt);
			\shade [ball color=black]  (1,0) circle (4pt);
			\draw (-1,0) -- (1,0) -- (0,2) -- cycle;
			\node[above] (z2) at (0,2) {$z_2$};
			\node[left] (x) at (-1,0) {$x$};
			\node[right] (y) at (1,0) {$y$};
			\node[below] (xy) at (0,0) {$\lambda-\lambda_1$};
			\node (xz2) at (-1.3,1) {$\lambda_1$};
			\node (yz2) at (1.3,1) {$\lambda_1$};
			\node (Gz1) at (0,-1.7) {$G[z_1]$};
			\end{tikzpicture}
		\end{tabular}
		\caption{\label{fig:Multiplicity}Multiplicity of edges in $G[x]$, $G[y]$ and $G[z_1]$.}
	\end{center}
	\vspace{-4ex}
\end{figure}
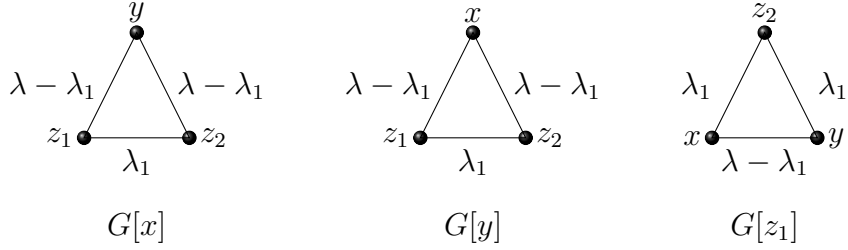

	As a consequence, in $G[z_1]$ the edges $xz_2$ and $yz_2$ appear with multiplicity $\lambda_1$ and the edge $xy$ appears with multiplicity $\lambda-\lambda_1$. 
    Considering the degree of $z_2$ in $G[z_1]$, we have that $2\lambda_1\leq\lambda$, and considering the degree of $y$ in $G[z_1]$ we have that $2(\lambda - \lambda_1) \leq \lambda$, which implies $\lambda \leq 2\lambda_1$. 
    Therefore $\lambda_1 = \lambda-\lambda_1 =\lambda/2$, and so each edge of $C$ appears with multiplicity $\lambda/2$.	
	We further note that this means that all of the pairs $\{xy, xz_1, xz_2,yz_1,yz_2,z_1z_2\}$ have now been covered $\lambda$ times and so cannot appear in blocks other than $B_1$, $B_2$, $B_3$, $B_4$. 
	
	Now suppose that Xeno plays $x$ as his first move, so Ophelia is forced to play $y$. Now, Xeno plays $z_1$, forcing Ophelia to play on $z_2$. At this point, Ophelia has played on $y$ and $z_2$, but we know that the edge $y\,z_2$ only appears in blocks containing either $x$ or $z_1$ ($B_2$, $B_4$), both of which Xeno has played on. Thus Ophelia is not threatening a win and Xeno may play where he pleases. 
	We suppose that Xeno plays some vertex $u \notin \{x,y,z_1,z_2\}$ (recalling that $v>5$).
	The blocks containing $\{u,x\}$ and $\{u,z_1\}$ cannot contain $y$ or $z_2$, since the edges $\{xy, xz_2, z_1y, z_1z_2\}$ only come from those block $B_i$, $1\leq i\leq 4$. Further, since $xz_1$ has already appeared $\lambda$ times, $\{u, x, z_1\}\notin\cB$.
	
	Each of the pairs $ux$ and $uz_1$ must appear $\lambda$ times in the blocks. So there must exist $a, b\in X$, $a\neq b$, such that $\{u,x,a\}$ and $\{u,z_1,b\}\in\cB$ (which implies $v>6$), but now Xeno has a scissor on $u$. 
	This contradicts the assumption that Ophelia can force a draw on $(X, \cB)$, and the result follows.
\end{proof}

We note that Theorem~\ref{TS(v,lambda)} shows that Xeno can win in the \ttt{} game on TS$(v,\lambda)$, and so Maker can win the corresponding \MB{} game by Theorem~\ref{MB to ttt}.

\subsection{Transversal Designs}
\label{TD}

We now investigate playing \ttt{} on a $\mathrm{TD}(k,n)$ for small values of $k$.

\begin{theorem}\label{TD(2,n)}
Xeno can win \ttt{} on a $\mathrm{TD}(2,n)$ if and only if $n\geq 2$.
\end{theorem}
\begin{proof}
	If $n=1$, there is only one block, which has cardinality $2$. Ophelia can force a draw on her turn by playing on that block. 
	Otherwise, $|G_1|=|G_2|=n\geq 2$, and Xeno plays on some point $x$ of the group $G_1$. Now, for every $y\in G_2$, there is a block $\{x,y\}$, Ophelia can block at most one such block, so since $n \geq 2$, Xeno wins the game on the next move. 
\end{proof}

Next, we show that any $\mathrm{TD}(3,n)$ is a Xeno win when $n \geq 3$. In the case where $n=3$ and the design is also resolvable, an RTD(3,3) can be thought of as the affine plane $\pi_3$ with a parallel class removed.
We note in passing that if we remove another parallel class from the RTD(3,3), it is equivalent to regular \ttt{} with the diagonals removed, and so it is not difficult to see that Ophelia can force a draw. 

\begin{theorem}
\label{k=3 is Xeno win}
	\label{TD(3,n)}
	Xeno can win \ttt{} on a $\mathrm{TD}(3,n)$ if and only if $n \geq 3$.
\end{theorem}
\begin{proof}
	For $n \leq 2$, Ophelia can force a draw by Theorem~\ref{thm:MB to ttt}. 
	If $n \geq 3$, we may assume that Xeno plays $X_1$ in group $G_1$. 
	The groups $G_2$ and $G_3$ are effectively equivalent for Ophelia's first turn, and so the course of play now splits depending on whether Ophelia plays in $G_1$ or not.
	We give a strategy for Xeno that ensures a win in either case.

	{\bf Case 1:} If Ophelia plays $O_1 \in G_1$, Xeno responds by playing in the same group again, $X_2\in G_1$ ($n\geq 3$). 
	If Ophelia plays on $O_2\in G_1$, then Xeno responds by playing $X_3\in G_2$. Note that $X_3$ is a scissor for Xeno, since
	the block containing $X_1X_3$ and the block containing $X_2X_3$ are both potential wins for Xeno; see Figure~\ref{fig: setup1}.
	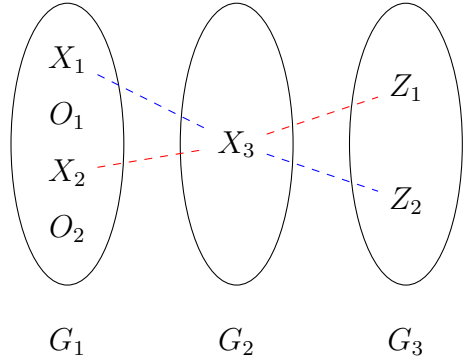
\begin{figure}[ht]
		\begin{center}
			\begin{tikzpicture}[scale=0.75, x = 1cm, y=1cm]
				\foreach \x in {1,4,7} {
					\draw (\x,4) ellipse (1 and 2.5);
				}
				\node (X1) at (1,5.5) {$X_1$};
				\node (O1) at (1,4.5) {$O_1$};
				\node (X2) at (1,3.5) {$X_2$};
				\node (O2) at (1,2.5) {$O_2$};
				\node (X3) at (4,4) {$X_3$};
				\node (Z1) at (7,5) {$Z_1$};
				\node (Z2) at (7,3) {$Z_2$};
				\node (G1) at (1,0.5) {$G_1$};
				\node (G2) at (4,0.5) {$G_2$};
				\node (G3) at (7,0.5) {$G_3$};
				\draw[color=red, dashed, postaction={decorate}] (X2) -- (X3) -- (Z1);
				\draw[color=blue, dashed, postaction={decorate}] (X1) -- (X3) -- (Z2);
			\end{tikzpicture}
		\caption{Playing \ttt{} on a $\mathrm{TD}(3,n)$, where Xeno and Ophelia played their first two moves on $G_1$, then Xeno plays his third in $G_2$.}\label{fig: setup1}
		\end{center}
	\end{figure}

	If Ophelia instead plays $O_2$ in $G_2$ (or equivalently, in $G_3$), Xeno is forced to play on the block containing $O_1O_2$ to avoid an Ophelia win, hence $X_3\in G_3$ (respectively $G_2$). But $X_3$ is a scissor for Xeno, and so Xeno wins on the next turn; see Figure~\ref{fig: setup2}.

	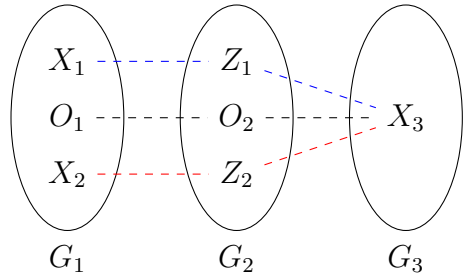
\begin{figure}[ht]
		\begin{center}
			\begin{tikzpicture}[scale=0.75, x = 1cm, y=1cm]
				\foreach \x in {1,4,7} {
					\draw (\x,4) ellipse (1 and 2);
				}
				\node (X1) at (1,5) {$X_1$};
				\node (O1) at (1,4) {$O_1$};
				\node (X2) at (1,3) {$X_2$};
				\node (O2) at (4,4) {$O_2$};
				\node (X3) at (7,4) {$X_3$};
				\node (Z1) at (4,5) {$Z_1$};
				\node (Z2) at (4,3) {$Z_2$};
				\node (G1) at (1,1.5) {$G_1$};
				\node (G2) at (4,1.5) {$G_2$};
				\node (G3) at (7,1.5) {$G_3$};
				\draw[color=blue, dashed, postaction={decorate}] (X1) -- (Z1) -- (X3);
				\draw[color=red, dashed, postaction={decorate}] (X2) -- (Z2) -- (X3);
				\draw[color=black, dashed, postaction={decorate}] (O1) -- (O2) -- (X3);
			\end{tikzpicture}
		\end{center}
		\caption{Playing \ttt{} on a $\mathrm{TD}(3,n)$, where Xeno played his first two moves in $G_1$ and Ophelia played her first move in $G_1$ and second move in $G_2$. Xeno blocks a win by playing in $G_3$, which is also a scissor for Xeno.}\label{fig: setup2}
	\end{figure}

	{\bf Case 2:} We now suppose that Ophelia does not play her first move in the same group as Xeno, so without loss of generality we may assume $O_1 \in G_2$. Xeno responds by playing in the other group that Ophelia did not play in, so $X_2 \in G_3$ with $\{X_1,O_1,X_2\} \notin \mathcal{B}$. 
	But now Ophelia must play her second move in the block containing $X_1X_2$, thus $O_2 \in G_2$ is forced. Now Xeno plays any remaining $X_3\in G_2$ ($n \geq 3$), any of which is a scissor for Xeno;
	see Figure~\ref{fig: setup3}.
	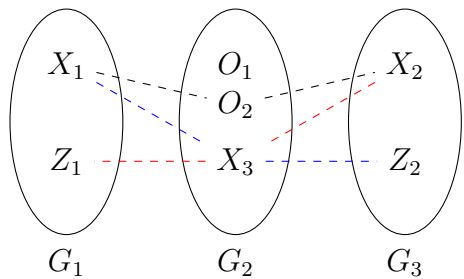
\begin{figure}
		\begin{center}
			\begin{tikzpicture}[scale=0.75, x = 1cm, y=1cm]
				\foreach \x in {1,4,7} {
					\draw (\x,4) ellipse (1 and 2);
				}
				\node (X1) at (1,5) {$X_1$};
				\node (O1) at (4,5) {$O_1$};
				\node (X2) at (7,5) {$X_2$};
				\node (O2) at (4,4.3) {$O_2$};
				\node (X3) at (4,3.3) {$X_3$};
				\node (Z1) at (1,3.3) {$Z_1$};
				\node (Z2) at (7,3.3) {$Z_2$};
				\node (G1) at (1,1.5) {$G_1$};
				\node (G2) at (4,1.5) {$G_2$};
				\node (G3) at (7,1.5) {$G_3$};
				\draw[dashed, postaction={decorate}] (X1) -- (O2);
				\draw[dashed, postaction={decorate}] (O2) -- (X2);
				\draw[color=red, dashed, postaction={decorate}] (X2) -- (X3) -- (Z1);
				\draw[color=blue, dashed, postaction={decorate}] (X1) -- (X3) -- (Z2);
			\end{tikzpicture}
		\end{center}
	\caption{Playing \ttt{} on a $\mathrm{TD}(3,n)$, where Ophelia plays her first and second move on a different group to Xeno's first move. }\label{fig: setup3}
	\end{figure}
\end{proof}

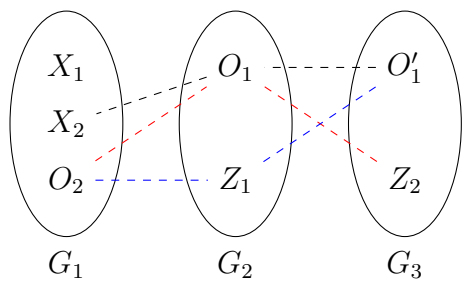
\begin{figure}
	\begin{center}
		\begin{tikzpicture}[scale=0.75, x = 1cm, y=1cm]
		\foreach \x in {1,4,7} {
			\draw (\x,4) ellipse (1 and 2);
		}
		\node (X1) at (1,5) {$X_1$};
		\node (O1) at (4,5) {$O_1$};
		\node (O2) at (7,5) {$O_1'$};
		\node (X2) at (1,4) {$X_2$};
		\node (Z1) at (4,3) {$Z_1$};
		\node (Z2) at (7,3) {$Z_2$};
		\node (O3) at (1,3) {$O_2$};
		\node (G1) at (1,1.5) {$G_1$};
		\node (G2) at (4,1.5) {$G_2$};
		\node (G3) at (7,1.5) {$G_3$};
		\draw[dashed, postaction={decorate}] (O2) -- (O1) -- (X2);
		\draw[color=red, dashed, postaction={decorate}] (O3) -- (O1) -- (Z2);
		\draw[color=blue, dashed, postaction={decorate}] (O3) -- (Z1) -- (O2);
		\end{tikzpicture}
	\end{center}
\caption{Playing \ttt{} on a $\mathrm{TD}(3,n)$, where Ophelia plays her second move on the same group to Xeno's first move, but not her first move.}\label{fig:extramove}
\end{figure}

A common question to ask is the following: On a Xeno win hypergraph, how close is Ophelia to winning the game? More precisely, how many extra moves would we need to give Ophelia to switch a Xeno win to an Ophelia win. 
We note that in a $\mathrm{TD}(3,n)$, if Ophelia gets one extra move on her first turn, she can win the game. 
To see this, suppose that
Xeno chooses $X_1\in G_1$ as his first move, Ophelia chooses two points, $O_1\in G_2$ and $O_1'\in G_3$, such that $\{X_1,O_1, O_1'\}$ is not a block of the design; see Figure~\ref{fig:extramove}. 
Now, Xeno is forced to play in $G_1$ to block an Ophelia win in the block containing $O_1$ and $O_1'$. Ophelia now plays in $O_2\in G_1$, which is a scissor for Ophelia, and she wins the game.

We note that Xeno's strategy given in the proof above is not 
a score-optimizing strategy. In fact, if he only plays the score-optimizing strategy, while stopping any winning moves from Ophelia, Ophelia can force a draw on \emph{atomic} $\mathrm{TD}(3,n)$s, which we now define.

Let $(X,\cG, \cB)$ be a $\mathrm{TD}(3, n)$ with $\cG =\{G_1, G_2, G_3\}$. 
Specify two points $x,y \in G_1$ and consider the graph $G[x] \cup G[y]$, recalling that $G[x]$ is the (multi-)graph with edge $ab$ for each $\{x,a,b\} \in \mathcal{B}$. 
As $x$ and $y$ both appear with each point of $G_2 \cup G_3$ exactly once each, the graph $G[x] \cup G[y]$ is $2$-regular. 
If $G[x] \cup G[y]$ consists of a single cycle, then the $\mathrm{TD}(3, n)$ is called \emph{Hamiltonian} with respect to $x$ and $y$. 

If every pair of points in the group $G_1$ generates a single cycle by the above process, then we say that the $\mathrm{TD}(3, n)$ is \emph{pan-Hamiltonian} with respect to $G_1$. If the $\mathrm{TD}(3, n)$ is pan-Hamiltonian with respect to each of its groups $G_1$, $G_2$, and $G_3$, then the $\mathrm{TD}(3, n)$ is said to be \emph{atomic}.
Atomic transversal designs have been studied in the context of Latin squares and perfect one-factorizations, see \cite{Rosa}. 
Results about atomic transversal designs imply that an atomic $\mathrm{TD}(3,n)$ exists for infinitely many values of $n$, see \cite{BryantMaenWan, WanPer, Wan2}. 

\begin{theorem}
	\label{pan-hamiltonian}
	Suppose we play \ttt{} on an atomic $\mathrm{TD}(3,n)$ with $n \geq 4$.
	If  Xeno plays a score-optimizing strategy,  then Ophelia can force a draw (by playing a score-optimizing strategy herself).
\end{theorem}
\begin{proof}
	Given an atomic $\mathrm{TD}(3,n)$ with groups $\cG =\{G_1, G_2, G_3\}$, we may arbitrarily take Xeno's first move to be $X_1\in G_1$. 
	Ophelia's score-optimizing move is to play any point in either $G_{2}$ or $G_{3}$, so without loss of generality we take $O_1 \in G_2$. 
	Now, Xeno's score-optimizing move is to choose a point in $G_2$ different from $O_{1}$, say $X_2 \in G_2$.
	We note that the two points $O_{1}$, $X_{2}$ are in the same group $G_{2}$, so $G[O_1]\cup G[X_1]$ is a Hamiltonian cycle.
	We see that Ophelia is under threat on the block containing the points $X_1$ and $X_2$, so she chooses her second move $O_2 \in G_3$ such that $\{X_{1}$, $X_{2}, O_{2} \}$ is a block.
	Now, Xeno is under threat on the block containing the points $O_1$ and $O_2$, so he chooses $X_3 \in G_1$ such that $\{ O_{1}$, $O_{2}, X_{3} \}$ is a block. 
	In general, we note that after Ophelia's $i^{\rm th}$ move, Xeno is under threat on the block containing the points $O_i$ and $O_{1}$, with $O_i \in G_3$ and $O_1 \in G_2$. 
	So he must choose to play each responding move in $G_1$ until $i=n$ to prevent Ophelia winning. 
	Similarly, after Xeno's $i^{\rm th}$ move, Ophelia is under threat on the block containing the points $X_i$ and $X_{2}$ with $X_i \in G_1$ and $X_2 \in G_2$, so she will choose each responding move $O_i$ in $G_3$ until $i=n$. 
	After the $(n+1)$th round, Xeno will have played all of the points in $G_1$ and Ophelia will have played all of the points in $G_3$, so no block remains for Xeno to be able to win.
	Thus the game ends in a draw. 
\end{proof}

Thus, comparing the results of Theorems~\ref{k=3 is Xeno win} and \ref{pan-hamiltonian} we can see that if Xeno plays a score-optimizing strategy, then Ophelia may be able to force a draw, depending on the structure of the design. 
This answers Question \ref{qn:score_opt} in the negative.

\section{Conclusion}

We have considered playing \ttt{} and \MB{} on BIBDs and transversal designs. We note that there is an interplay between the games by Theorem~\ref{MB to ttt}. 
Corollary~\ref{TD 2^k > n(n+1)} shows that when playing \MB{} on a TD$(k,n)$, Maker can win when $n > k2^{k-3}$ and Breaker can win when $n \leq 2^{k/2} - 1/2$, which means that Ophelia can force a draw in \ttt{} when $n \leq 2^{k/2} - 1/2$. 
For general BIBD$(v,k,\lambda)$, we have shown that Maker can win \MB{} when $v >k(k-1)2^{k-3}+1$ by Lemma~\ref{lem:BIBD MB}. 
This lemma also shows that Breaker can win when $v< \frac{-k+1 + \sqrt{(k+1)^2 + k(k-1)2^{k+2}}}{2}$, and hence Ophelia can force a draw in this case.
Perhaps most importantly, we have shown in Theorem~\ref{TS(v,lambda)} that each TS$(v, \lambda)$ is a Xeno win in \ttt{} if and only if $v\geq 5$.

Particular results are summarized for \MB{} in Table~\ref{Maker_Breaker} and those for \ttt{} in Table~\ref{Tic-Tac-Toe}.
However, we note that in order for a $\mathrm{TD}(k,n)$ to exist, $k\leq n+1$, and while existence is known for all prime power orders of $n$, there are many pairs ($k,n)$ with $k\leq n+1$ for which existence of a $\mathrm{TD}(k,n)$ is unknown. Similarly, for a BIBD$(v, k,\lambda)$ to exist, the necessary conditions of Equation~(\ref{bibd nec}) from Subsection~\ref{subs Designs}
must be satisfied. In particular, when $\lambda=1$, this means that $v\equiv 1$ or $k\bmod{k(k-1)}$. 

\renewcommand{\arraystretch}{1}

\begin{table}[ht]
	\begin{center}
		\begin{tabular}{|l|l|l|l|} \hline
			Designs & Maker & Breaker & Reference\\ \hline\hline
			TD$(k,n)$&$ n > k2^{k-3}$&$n \leq 2^{k/2} - 1/2$ & Corollary~\ref{TD 2^k > n(n+1)} \\ \hline
			TD$(2, n)$ &$n \geq 2$ &$n=1$ &Corollary~\ref{cor:final_TD34} \\ \hline
			TD$(3, n)$&$n \geq 3$ & $n\leq2$& Corollary~\ref{cor:final_TD34} \\ \hline
			 TD$(4, n)$&$n \geq9$&$n=3,4$ & Corollary~\ref{cor:small values of k} \\
			 & & & Theorem~\ref{TD(4,4)}\\ \hline
			 TD(4,5)&$\;\;\;\;\checkmark$&&Computer search \\ \hline \hline
			BIBD$(v,3,1)$ & $v \geq7$ &$v=3$ &Corollary~\ref{cor:BIBD(v,3,1} \\ \hline
			BIBD$(v,4,1)$ & $v\geq16$ & $v=4, 13$& Theorem~\ref{BIBD(v,4,1)} \\
			 \hline \hline
			TS$(v,\lambda)$ & $v\geq 5$ &$v=3,4$ & Corollary~\ref{cor:TS(v,lambda)-MB} \\ \hline
		\end{tabular}
		\caption{Results for the \MB{} Game. \label{Maker_Breaker}}
	\end{center}
\end{table}

\begin{table}[ht]
	\begin{center}
		\begin{tabular}{|l|l|l|l|} \hline
			Designs& Xeno &Ophelia & Reference\\ \hline\hline
			TD$(k,n)$&&$n \leq 2^{k/2} - 1/2$ & Theorem~\ref{thm:MB to ttt} \\ \hline
			TD$(2, n)$ & $n \geq 2$ & $n = 1$&Theorem~\ref{TD(2,n)} \\ \hline
			TD$(3, n)$&$n\geq3$&$n \leq 2$& Theorem~\ref{k=3 is Xeno win} \\ \hline
			TD$(4, n)$&&$n=3,4$&Theorem~\ref{thm:MB to ttt} \\ \hline
			TD$(4, 5)$ & & $\;\;\;\;\checkmark$ & Computer search \\ \hline \hline
			BIBD$(v,3,1)$ & $v\geq 7$ & $v=3$ & Theorem~\ref{STS(v)}\\ \hline
			BIBD$(v,4,1)$ & & $v=4, 13$& Theorem~\ref{thm:MB to ttt}\\ \hline\hline 
			TS$(v,\lambda)$&$v\geq5$&$v=3,4$& Theorem~\ref{TS(v,lambda)} \\ \hline
		\end{tabular}
		\caption{Results for the \ttt{} Game. \label{Tic-Tac-Toe}}
	\end{center}
\end{table}

Theorem~\ref{thm:MB to ttt} summarizes the consequences of our results in Section \ref{Weak Breaker win} on \MB{} to \ttt{}. 
Given the results of Section~\ref{Strong}, we are able to give the following results on \MB{}.
\begin{corollary} \label{cor:final_TD34}
     For $k\in\{2,3\}$, playing \MB{} on a $\mathrm{TD}(k,n)$ is Maker win if and only if $n \geq k$.
\end{corollary}
\begin{proof}
    This is found by applying Theorem~\ref{MB to ttt} to Theorems~\ref{TD(2,n)} and \ref{TD(3,n)}. 
\end{proof}
\begin{corollary}\label{cor:TS(v,lambda)-MB}
    Playing \MB{} on a {\rm TS}$(v,\lambda)$ is Maker win if and only if $v\geq 5$.
\end{corollary}
\begin{proof}
    The cases $v\geq 5$ follow by applying Theorem~\ref{MB to ttt} to Theorem~\ref{TS(v,lambda)}. The cases  $v=3,4$ are obvious Breaker wins as Maker has fewer than three moves and so cannot complete a block.
\end{proof}

In Section \ref{Strong}, we analyzed the \ttt{} game on several designs, including Transversal designs TD$(k,n)$, in particular for $k=2,3$. 
When $k=4$, very little is known. We know that the $\mathrm{TD}(4,3)$ and the $\mathrm{TD}(4,4)$ are Ophelia draw, which both follow from Theorem~\ref{thm:MB to ttt}.
We have verified by computer search that the $\mathrm{TD}(4,5)$ is an Ophelia draw. 
Very surprisingly, this computation also found that the $\mathrm{TD}(4,5)$ is Maker win, which yields the first known example of a hypergraph that is Maker win where Ophelia can force a draw. 
In this hypergraph, we have verified that the removal of any vertex breaks this property, but there are three hyperedges that can be removed while retaining the property. We conjecture that this is the smallest order for a graph with this property.
\begin{conjecture}
    Any hypergraph that is Maker win and where Ophelia can force a draw must have at least $20$ vertices. 
\end{conjecture}

Currently, we have no examples of a $\mathrm{TD}(k,n)$ with $k>3$ where Xeno wins with the \ttt{} game, and the existence of such Xeno-win $\mathrm{TD}(k,n)$ remains an open question. 

We note that Theorem~\ref{TD(4,4)} with Theorem~\ref{MB to ttt} implies that when playing \ttt{} on a \TD$(4,4)$, Ophelia can force a draw.
As mentioned, $\pi_4$ can be obtained by taking by adding the set of groups of the TD$(4,4)$ to the set of blocks.
We know that $\pi_4$ is a Xeno win by Theorem~\ref{Affine and Projective Planes Result}. 
That is, $\pi_4$ is Xeno-win for \ttt{} and hence Maker win for \MB{}, however a \TD$(4,4)$ is Breaker win for \MB{} and hence is Ophelia-draw for \ttt{}. 
So in some sense, these cases lie at the boundary between Xeno win and Ophelia draw. This boundary has been called the extra edge paradox, and other (uniform) examples can be found in \cite{MR2706074}. 
Up to isomorphism, there is only one \TD$(4,4)$ and only one $\pi_4$. We have verified using a computer that:
\begin{itemize}
	\item
	adding one group as a block is still a Breaker win, hence also an Ophelia draw in \ttt{};
	\item 
	adding two groups as blocks is still a Breaker win, hence also an Ophelia draw in \ttt{};
	\item 
	adding three groups as blocks is a Xeno win in \ttt{}, hence also a Maker win in \MB{}.
\end{itemize}

We can thus see that designs seem to be a good place to look for extreme cases. 
While we have determined the outcomes for many games, there is still much work to be done.

\bibliographystyle{siam}
\bibliography{TTT}

\end{document}